\documentclass[3p,11pt]{elsarticle}
\usepackage{amsmath}
\usepackage{appendix}
\usepackage{stmaryrd}

\usepackage{amsfonts}
\usepackage{amssymb}
\usepackage{mathrsfs}
\usepackage{amsfonts}
\usepackage{graphicx}
\usepackage{color}
\usepackage{graphics,color,epsfig}
\usepackage{stmaryrd}
\usepackage{enumitem}
\usepackage{amsthm}
\usepackage{subcaption}
\usepackage{caption}
\usepackage{amsthm}
\usepackage{chngcntr}
\usepackage{hyperref}

\usepackage{multirow}

\usepackage{booktabs}
\usepackage{xcolor} 
\usepackage{graphicx} 
\usepackage{siunitx}

\usepackage{tabularx}       
\usepackage{siunitx}        
\usepackage{ragged2e}       

\usepackage{float}
\usepackage{makecell}

 \usepackage[linesnumbered,ruled,vlined]{algorithm2e} 
\usepackage{bm} 

\newtheorem{cor}{Corollary}[section]

\usepackage{bm}

\newcommand{\dx}[1][x]{\,{\rm d}#1}

\newtheorem{theorem}{Theorem}[section]
\newtheorem{lemma}{Lemma}[section]
\newtheorem{example}{Example}[section]
\newtheorem{remark}{Remark}[section]
\newtheorem{defn}{Definition}[section]
\newcommand \dint {\displaystyle\int}
\newcommand{\bs}[1]{\boldsymbol{#1}}
\def\B{\mathbb{B}}
\def \bx{\bs x}
\def \by{\bs y}
\def \bz{\bs z}
\def \bY{\bs Y}
\def\d{\mathrm{d}}

\def \Q{\widetilde{Q}}

\numberwithin{equation}{section} 
\usepackage{chngcntr} 
\counterwithin{algocf}{section} 

\journal{Journal of Computational Physics}
 \graphicspath{{figures/}}
 
\begin{document}
\begin{frontmatter}

\title{FNWoS: Fractional Neural Walk-on-Spheres Methods for High-Dimensional  PDEs Driven by $\alpha$-stable L\'{e}vy Process on Irregular Domains}
\author[SNH]{Ling Guo}
\author[SDU]{Mingxin Qin\corref{cor}}
\cortext[cor]{Corresponding Author}\ead{mingxin\_qin@163.com}
\author[SUFE]{Changtao Sheng}
\author[SJTU]{Hao Wu}
\author[SDU]{Fanhai Zeng}
\address[SNH]{Department of Mathematics, Shanghai Normal University, Shanghai, China}
\address[SDU]{School of Mathematics, Shandong University, Jinan, China}
\address[SUFE]{School of Mathematics, Shanghai University of Finance and Economics, Shanghai, China}
\address[SJTU]{School of Mathematical Sciences, Institute of Natural Sciences, and MOE-LSC, Shanghai Jiao Tong University, Shanghai, China}

\begin{abstract}
In this paper, we develop a highly parallel and derivative-free fractional neural walk-on-spheres method (FNWoS) for solving high-dimensional fractional Poisson equations on irregular domains. We first propose a simplified fractional walk-on-spheres (FWoS) scheme that replaces the high-dimensional normalized weight integral with a constant weight and adopts a correspondingly simpler sampling density, substantially reducing per-trajectory cost. To mitigate the slow convergence of standard Monte Carlo sampling, FNWoS is then proposed via integrating this simplified FWoS estimator, derived from the Feynman-Kac representation, with a neural network surrogate. By amortizing  sampling effort over the entire domain during training, FNWoS  achieves more accurate  evaluation at arbitrary query points with dramatically fewer trajectories than classical FWoS. To further enhance efficiency in regimes where the fractional order $\alpha$ is close to 2 and trajectories become excessively long, we introduce a truncated path strategy with a prescribed maximum step count. Building on this, we propose a buffered supervision mechanism that caches training pairs and progressively refines their Monte Carlo targets during training, removing the need to precompute a highly accurate training set and yielding the buffered fractional neural walk-on-spheres method (BFNWoS). Extensive numerical experiments, including tests on irregular domains and problems with dimensions up to $1000$, demonstrate the accuracy, scalability, and computational efficiency of the proposed methods.
\end{abstract}


\end{frontmatter}
 
 \section{Introduction}
 
 Anomalous diffusion arises in a wide range of natural systems, including viscoelastic materials, fractals, dispersion in porous media, and biological structures (see, e.g.,\cite{Benson2006,MR3413590,Kusnezov1999,Feder1996} and the references therein). The fractional Laplacian operator $(-\Delta)^{\frac{\alpha}{2}}$ serves as a powerful analytical framework for capturing such behaviors, thereby motivating extensive research on fractional PDEs for modeling purposes. In direct analogy with the classical Laplacian, which generates Brownian motion, the fractional Laplacian is the infinitesimal generator of a standard isotropic $\alpha$-stable L\'{e}vy process.
 In this paper, we consider the following high-dimensional fractional Poisson equation with global Dirichlet boundary conditions on irregular domains  
\begin{equation}\label{fractional-laplace}
    \begin{aligned}
        \begin{cases}
        (-\Delta)^{\frac{\alpha}{2}} u(\bx) = f(\bx), & \quad \text{in } \Omega, \\[2pt]
        u(\bx) = g(\bx), & \quad \text{in } \mathbb{R}^d \setminus \Omega,
        \end{cases}
    \end{aligned}
\end{equation}
where $\alpha \in (0,2)$, $\Omega$ is an open bounded domain, and the fractional Laplacian is defined as (cf.\,\cite{di2012hitchhiker})
\begin{equation}
(-\Delta)^{\frac{\alpha}{2}}u(\bx) = C_{d,\alpha} \, \text{p.v.} \int_{\mathbb{R}^d} \frac{u(\bx) - u(\by)}{|\bx-\by|^{d+\alpha}} \, \d \by, \quad C_{d,\alpha} := \frac{\alpha2^{\alpha-1} \Gamma\left(\frac{d+\alpha}{2}\right)}{\pi^{\frac{d}{2}} \Gamma\left(1-\frac{\alpha}{2}\right)},
\end{equation}
where “p.v.” stands for the Cauchy principal value of the integral. 

The numerical solution of the fractional Poisson equation \eqref{fractional-laplace} is challenging due to the nonlocality of the fractional Laplacian, the hypersingular nature of the integrals, and the low regularity of the solution near the boundary. Substantial efforts have been devoted to developing numerical algorithms for fractional PDEs. The first category comprises traditional deterministic methods, including finite difference schemes~\cite{hao2021fractional,minden2020simple,duo2019accurate,li2012finite,huang2014numerical}, finite element approaches~\cite{acosta2017fractional,bonito2019numerical,ainsworth2018towards,ainsworth2017aspects,bu2014galerkin}, and spectral techniques~\cite{sheng2020fast,song2017computing,hao2021sharp,tang2020rational,xu2018spectral}. Although these methods perform well on simple bounded domains, their efficiency deteriorates markedly when applied to complex geometries, even in two or three dimensions. To address this issue, Antil et al.\,\cite{antil2021approximation} proposed solving the fractional Laplacian on hypercubes in the frequency domain by exploiting the explicit Fourier transform of the sinc function, and further combined it with the fictitious domain technique to enable computations via FFT on irregular domains in 2D and 3D. Huang et al.\,\cite{huang2024} introduced the grid-overlay finite-difference method combined with an FFT-based solver to handle fractional PDEs on irregular domains in 2D and 3D. However, deterministic approaches still face an insurmountable challenge posed by the curse of dimensionality when extended to higher-dimensional settings (i.e. for $d>3$).

Unlike deterministic approaches, machine learning provides an effective means to mitigate the curse of dimensionality due to its superior capability in representing high-dimensional functions (see, e.g., \cite{e2017deep,raissi2018Forward,han2020solving,Cai2025Soc,Cai2025Martingale} and the references therein). A representative example is the physics-informed neural networks (PINNs) \cite{raissi2019physics}, which is well suited for such problems as it can efficiently handle classical PDEs by leveraging mature tools such as automatic differentiation. However, the fractional Laplacian operators cannot directly support automatic differentiation due to their intrinsic nonlocal nature. To address the inapplicability of automatic differentiation for PDEs with nonlocal operators (such as fractional PDEs, FPDEs), Pang et al. \cite{pang2019fpinns} extended PINNs to fractional PINNs (fPINNs) for solving FPDEs. Subsequently, Ma et al. \cite{ma2023bi} proposed the biorthogonal fPINNs to tackle stochastic fractional PDEs. Nevertheless, the discretization employed in fPINNs, which relies on numerical schemes such as finite differences, remains susceptible to the curse of dimensionality.
To alleviate the curse of dimensionality in high-dimensional fractional PDEs, Guo et al. \cite{guo2022monte} exploited the nonlocal property of fractional operators to reformulate them as integral problems and constructed efficient sampling strategies based on the associated density functions, combining them with Monte Carlo methods to approximate the fractional derivatives of deep neural network outputs. However, approximating integrals with singularities in fractional derivatives via Monte Carlo introduces inherent variance and potential errors, which hinder the convergence of PINNs. To overcome this issue, Hu et al. \cite{hu2024tackling} employed Gauss quadrature to evaluate the integrals, providing a theoretical guarantee of accuracy and yielding a more precise algorithm.

Another natural approach to addressing high-dimensional fractional problems is to employ stochastic algorithms, which provide a distinctive framework for the numerical solution of high-dimensional PDEs, with significant progress reported in recent studies \cite{Shao2020,Lei2025}. Among various stochastic approaches, those based on the Feynman-Kac representation are particularly representative. Their core idea is to use the Feynman-Kac formula to establish a probabilistic connection between PDEs and stochastic processes, while incorporating Monte Carlo techniques to effectively mitigate the curse of dimensionality. In this context, under mild regularity assumptions on the data, namely that \(g(x)\) is continuous with \(g \in L^{1}_\alpha(\Omega^c)\), and \(f \in C^{\alpha+\varepsilon}(\overline{\Omega})\) for some \(\varepsilon>0\), the fractional Poisson problem~\eqref{fractional-laplace} admits a unique continuous solution, which can be expressed via the Feynman-Kac formula (see~\cite[Thm.\,6.1]{Kyprianou18}):
\begin{equation}\label{Feynmankac}
u(\bx) = \mathbb{E}_{X_{0}^{\alpha }=\bx}\big[g(X^{\alpha}_{\tau_{\Omega}})\big] 
        + \mathbb{E}_{X_{0}^{\alpha }=\bx}\Big[\int_{0}^{\tau_{\Omega}} f(X^{\alpha}_{s})\,{\rm d}s\Big],
      \end{equation}
where $\tau_{\Omega} = \inf \{ t>0 : X_{t}^{\alpha}\notin \Omega \}$ and $\{X_{t}^{\alpha}\}_{t\geq0}$ is a symmetric $\alpha$-stable L\'{e}vy process with $X_{0}^{\alpha}=\bx \in \Omega$. In recent years, there has been growing interest in stochastic simulation approaches for nonlocal equations. In general, methods based on the Feynman-Kac representation can be grouped into two major directions. The first is to reformulate the original equation into a system of backward stochastic differential equations; for instance, Yang et\,al.~\cite{Yang2022a} proposed a stochastic scheme for semilinear nonlocal diffusion equations with volume constraints. The second direction adopts a more direct perspective, where numerical solutions are obtained by simulating the underlying stochastic process. A typical example is the work of Kyprianou et\,al.~\cite{Kyprianou18}, who proposed an efficient Monte Carlo algorithm for the fractional Poisson equation in high dimensions using the walk-on-spheres strategy derived from the Feynman-Kac formulation. In a related study, Jiao et al.~\cite{Jiao2022a} extended this strategy to high dimensions by introducing a quadrature rule for evaluating the integral representation over the ball and employing a rejection sampling method. More recently, Sheng et al.\,\cite{Sheng2023,Sheng2026} presented an efficient solution of fractional PDEs defined on irregular domains in high dimensions, accompanied by a rigorous numerical analysis that substantiates both its accuracy and computational efficiency. However, these stochastic algorithms generally achieve at most half-order accuracy in the sampling direction, so obtaining solutions of sufficient accuracy often entails a large number of samples. Hence, reducing the sampling cost while preserving the desired accuracy represents a central avenue for enhancing the efficiency of existing stochastic methods.

The main purpose of this paper is to develop a FNWoS method for solving PDEs driven by $\alpha$-stable Lévy processes on irregular domains in high dimensions. The core idea of the proposed method is to embed neural networks into the walk-on-spheres framework, exploiting stochastic path information to remove the reliance on automatic differentiation. More specifically, we train the DNN surrogate $v_\theta$ by minimizing
the following empirical loss function  
\begin{equation}\label{loss_introduc}
    \mathcal{L}_{\text{FNWoS}}[v_\theta] := \frac{1}{M}\sum_{j=1}^M \!\Big[ v_\theta(\bx_j) - u_\text{FWoS}(\bx_j) \Big]^2, 
\end{equation}
where $u_\text{FWoS}(\bx_j)$ denotes the simplified FWoS estimate,
based on stochastic path information, evaluated at the input points
$\{\bx_j\}_{j=1}^M$ (cf. Section~\ref{sec3}).
 Compared with traditional deterministic numerical methods such as the finite element method, finite difference method, and spectral method, the proposed approach is more suitable for solving fractional Poisson equations on complex geometries as well as in high-dimensional settings. Moreover, in comparison with existing meshfree methods or deep learning, the main innovations and contributions of this work are summarized as follows.
\begin{itemize}
    \item We  simplify the classical FWoS method proposed in \cite[Theorem~2.1]{Sheng2023} in two key aspects: {\rm(i)} the normalized weight $\omega_{r_k}$ is reduced from the integral form $\int_{\mathbb{B}_{r_k}^d} G_{r_k}(\bx,\by)\,\mathrm{d}\by$ to a constant $\omega_{r_k}$, thereby avoiding an expensive evaluation that becomes particularly burdensome in high-dimensional settings; {\rm(ii)} the sampling density $Q_{r_k}$ for generating the random variable $\{Y_{k+1}\}$ is replaced by a simpler distribution that is easier to evaluate and sample from. These modifications substantially reduce the computational cost of FWoS and enable more efficient generation of the required random variables.
    
    \vspace{-6pt}
    \item {We propose the FNWoS method, an extension of the simplified FWoS approach, which integrates neural networks to learn a surrogate $v_\theta$ for the solution of fractional Poisson equations. The resulting solver is derivative-free and fully parallelizable across the domain $\Omega$. By amortizing the computational cost of the simplified FWoS estimator during the training phase, FNWoS enables rapid evaluation of the solution at arbitrary query locations. For instance, as shown in Figure~\ref{fig:WoS_NWoS_ball_compare and fig:max_step}(b) for the $10$-dimensional unit ball, FNWoS achieves higher accuracy with only $100$ trajectories, whereas the simplified FWoS requires $10000$ trajectories trajectories to reach a comparable level of accuracy.}
    
    \vspace{-6pt}
    \item {To further improve efficiency, especially when $\alpha$ approaches $2$, we cap the maximum number of walking steps per trajectory to prevent excessive cost caused by overly long paths. Building on this truncation, we introduce a buffered supervision strategy that eliminates the need to generate a high-fidelity training set upfront and instead progressively refines the supervision signals during training. We refer to the resulting method as BFNWoS, and its efficiency gain is illustrated in  Figure~\ref{fig:time and rl2}.}
    
\end{itemize}

The rest of this paper is organized as follows. In Section~\ref{sec2}, we present a simplified probabilistic representation of the fractional Poisson equation and give a simplified FWoS scheme for high-dimensional problems on irregular domains. Section~\ref{sec3} develops the FNWoS method and its optimized variant, i.e., the BFNWoS method, two extensions of simplified FWoS that employ neural networks to approximate the solution of the fractional Poisson equation on high-dimensional irregular domains. Section~\ref{sec4} reports numerical experiments that demonstrate the efficiency and accuracy of the proposed approaches. Finally, concluding remarks are given in Section~\ref{sec5}.

 \section{Simplified FWoS scheme based on Feynman-Kac formula}\label{sec2}
In this section, we first introduce a simplified probabilistic representation of the fractional Poisson equation involving the integral fractional Laplacian, which simplifies the sampling process by separating the singular and smooth components of the Green’s function and serves as a fundamental building block of our algorithm. Based on the reformulated Feynman-Kac representation, we develop an efficient Monte Carlo method combined with the WoS method to solve the fractional Poisson equation on irregular domains in arbitrary high dimensions.

 \subsection{\bf{Feynman-Kac representation in expectation form.}}
We begin with the following explicit expressions for the Green’s function and the Poisson kernel on an arbitrary dimensional ball  associated with the $\alpha$-stable Lévy process (see \cite{MR3461641,Sheng2023}).
\begin{defn}\label{ddGHPF} {\bf (Green function and Poisson kernel).} 
Let $r>0$, for any $\bx,\by\in\B^{d}_{r}$ and $\bz\in\mathbb{R}^{d}\!\setminus \! \overline{\B}^{d}_{r}$, we define the Green function for $\bx\neq \by$ as
\begin{equation}\label{Gr}
    G_{r}(\bx,\by) = \begin{cases} \widetilde{C}^\alpha_d |\by-\bx|^{\alpha-d}\dint_{0}^{\varrho(\bx,\by)}\frac{t^{\alpha/2-1}}{(t+1)^{\frac{d}{2}}}\,\d t,\;\;&\alpha\neq d,\\[11pt]
   \widetilde{C}_0\log \Big{(}   \frac{ r^{2}-\bx \by+\sqrt{ (r^{2} -\bx^{2})(r^{2}-\by^{2}) }}{r|\by-\bx|}  \Big{)},\;\;&\alpha=d,
    \end{cases}
\end{equation}
and the Poisson kernel as
 \begin{equation}\label{Pr}
        P_r(\bx,\bz)=\widehat{C}^\alpha_d\Big(\frac{r^2-|\bx|^2}{|\bz|^2-r^2}\Big)^{\alpha/2}\frac{1}{|\bx-\bz|^d},\;\;\bx\in \B^{d}_{r},\;\;\bz\in \mathbb{R}^{d}\!\setminus \! \overline{\B}^{d}_{r},
    \end{equation} 
where 
\begin{equation}\begin{split}\label{cGreen}
 & \widetilde{C}^\alpha_d=\frac{\Gamma(d/2)}{2^{\alpha}\pi^{\frac d2}\Gamma^2(\alpha/2)},\qquad
 \widetilde{C}_0=\frac{\Gamma(1/2)}{2^{\frac{1}{2}}\pi^{\frac{1}{2}}\Gamma^2(1/4)},\\
 & \widehat{C}^\alpha_d=\frac{\Gamma(d/2)\sin(\pi\alpha/2)}{\pi^{\frac{d}{2}+1}},\;\;\;\varrho(\bx,\by) = \frac{(r^{2} - |\bx|^{2})(r^{2}-|\by|^{2})}{r^{2}|\bx-\by|^{2}}.
\end{split}\end{equation}
\end{defn}

We then derive an explicit expression for the solution of \eqref{fractional-laplace} on a ball of arbitrary dimensions, which is expressed in the form of an expectation.
\begin{lemma}\label{lem2.1}
Let $\Omega=\mathbb{B}_r^d$ with $r>0$ centered at the origin, and assume that $f \in C^{\alpha+\varepsilon}(\mathbb{B} _r^d) \cap C(\overline{\mathbb{B} _r^d})$ and $g \in L_\alpha ^1( \mathbb{R}^d\backslash \mathbb{B} _r^d)$; then 
the expectation form of the solution of \eqref{fractional-laplace} at the center $\bx$ of the ball can be expressed as
 \begin{equation}\label{eq-lem-2-2-1}
    u(\bx) = \mathbb{E}_{P_r}\big[g(Z)\big] + \omega_r \mathbb{E}_{{Q}_r}\Big[f(Y)W(\bx,Y)\Big], \quad Y \in \mathbb{B}_r^d, \, Z \in \mathbb{R}^d \setminus \mathbb{B}_r^d,
\end{equation}
where  $P_r$ is defined by \eqref{Pr}, the normalized weight function $\omega_r:=\omega_{r}(d,\alpha)= \frac{r^\alpha B(\frac{d-\alpha}{2},\frac{\alpha}{2})}{\alpha 2^{\alpha-1} \Gamma^2(\alpha/2)}$, and 
 \begin{equation}\label{omegaweight}
    \begin{aligned}
  Q_{r}(\bx,\by)= \frac{\Gamma(d/2)r^\alpha}{2\alpha\pi^{d/2}} |\boldsymbol{y} - \boldsymbol{x}|^{\alpha -d},\;\;\; 
      W(\bx,\by)=1-I\Big(\varrho^*(\boldsymbol{x},\boldsymbol{y});\frac{d-\alpha}{2},\frac{\alpha}{2}\Big),
    \end{aligned}
 \end{equation}
 with   $I(\cdot\,;\,\cdot,\cdot)$ denotes the incomplete Beta function, $\varrho(\bx,\by)$ is defined in \eqref{cGreen} and 
 $$\varrho^*(\boldsymbol{x},\boldsymbol{y})=\frac{r^2|\boldsymbol{x-y}|^2}{(r^2-|\boldsymbol{x}|^2)(r^2-|\boldsymbol{y}|^2)+r^2|\boldsymbol{x-y}|^2}.
 $$
\end{lemma}
 
\begin{figure}[htbp]
\centering 
\includegraphics[width=0.45\textwidth]{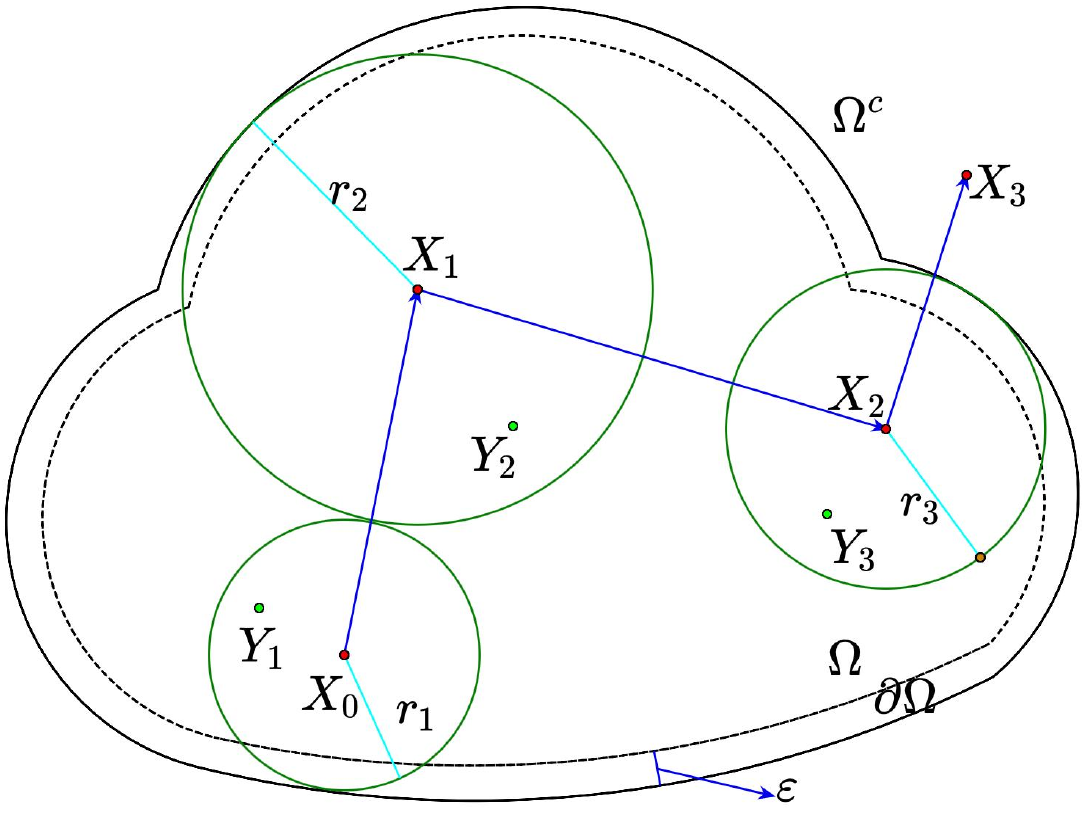} 
\caption{Random walk paths on a 2D irregular domain generated by the simplified FWoS.}
\label{WoS}
\end{figure}

We sketch the derivation of Lemma \ref{lem2.1} in Appendix A to avoid distraction from the main topic.  
Clearly, the above result is only applicable to arbitrary dimensional balls, as explicit expressions for the Green’s function and Poisson kernel are generally unavailable for irregular domains. To extend this approach to irregular domains, we incorporate the idea of the walk-on-spheres method, where a sequence of inscribed small balls is used to approximate the trajectory within the irregular domains (cf. Figure\,\ref{WoS}). Then, Lemma \ref{lem2.1} for arbitrary dimensional balls play a crucial role in computations within these small balls. Since the Markov process exits the given region in finite time, we define the stopping step for the random walk $X_\ell$ as  
$m^\ast= \inf\{\ell : X_\ell \notin \Omega\}.$ 
We present in Appendix C the detailed derivation of the following main results on the probabilistic representation of solutions over irregular domains.

\begin{theorem}\label{th-2-1}
Let $\alpha \in (0, 2)$, and $\Omega$ be an open bounded domain. Assume that $f \in C^{\alpha+\varepsilon}(\Omega) \cap C(\overline{\Omega})$, $g \in L_\alpha^1(\Omega^c) \cap C(\Omega^c)$, then the solution of \eqref{fractional-laplace} in $L_\alpha^1(\Omega)$ can be expressed as the form of the following conditional expectation
\begin{equation} \label{th-2-1-1}\begin{aligned}
    u(\bx) =&\;\,  \mathbb{E}_{P_{r_{m^*-1}}}\big[g(X_{m^*})\big]+\mathbb{E}_{{Q}_{r_k}}\bigg[\sum_{k=0}^{m^*-1} \omega_{r_k} f(Y_{k+1})W(\bx_k,Y_{k+1}) \Big| X_{\tau_k}^\alpha = \bx_k\bigg],
\end{aligned}\end{equation}
for any $Y_{k+1} \in \mathbb{B}_{r_k}^d, \, X_{m^*} \in \Omega^c.$
\end{theorem}

\begin{remark}
It is worth noting that the probabilistic representation of the solution presented above simplifies the result in {\rm\cite[Theorem 2.1]{Sheng2023}} in two aspects: {\rm(i)} the normalized weight function 
$\omega_{r_k}$ is simplified from the general form  $\int_{\mathbb{B}_{r_k}^d} G_{r_k}(\bx, \by)\, \mathrm{d}\by$ to a constant $\omega_{r_k}(d,\alpha)=\frac{{r_k}^\alpha B(\frac{d-\alpha}{2},\frac{\alpha}{2})}{\alpha 2^{\alpha-1} \Gamma^2(\alpha/2)}$ that depends on the radius $r_k$, where the  formulation
in {\rm\cite[Theorem 2.1]{Sheng2023}} is computationally expensive, especially in high-dimensional settings; 
{\rm(ii)} the density function  ${Q}_{r_k}$ associated with the random variable $\{Y_{k+1}\}$ is  also simplified.
\end{remark}

\subsection{\bf{Implementation of the simplified FWoS scheme}} For the simplified FWoS method, we only need to generate the locations of the stochastic process,
$\{X_k\}_{k=0}^{m^\ast}$, along with their random pairs $\{Y_k\}_{k=1}^{m^\ast}$, 
and substitute them into \eqref{th-2-1-1} to obtain the numerical solution. This implies that for a given initial position 
$\bx=(x_1,\cdots,x_d)$, using cartesian coordinates in high-dimensional settings would require a large number of integrals, which leads to excessively cumbersome expressions. Instead, by adopting spherical coordinates, the next position can be determined by computing 
the jump distance, while the direction on the sphere remains uniformly distributed. 
We now recall  the $d$-dimensional spherical coordinates 

\begin{equation*}
\begin{aligned}
x_1 &= \rho \cos \theta_1;~x_2 = \rho \sin \theta_1 \cos \theta_2;~\dots~;
~x_{d-1}= \rho \sin \theta_1 \cdots \sin \theta_{d-2} \cos \theta_{d-1};\\
x_d &= \rho \sin \theta_1 \cdots \sin \theta_{d-2} \sin \theta_{d-1},\;\;\;
\theta_1, \ldots, \theta_{d-2}\in [0, \pi],~\theta_{d-1}\in [0, 2\pi],
\end{aligned}
\end{equation*}
with  the spherical volume element 
\begin{equation*}
\d\bx = \rho^{d-1} \sin^{d-2} \theta_1 \sin^{d-3} \theta_2 \cdots \sin \theta_{d-2} \, \d\rho \, \d\theta_1 \, \d\theta_2 \cdots \d\theta_{d-1}.
\end{equation*}

In practical computations, to reduce the computational cost, we introduce a boundary layer of thickness $\varepsilon$ inside the domain. Once the small ball intersects this thickened boundary, it is treated as having exited the domain, thereby reducing the number of stopping steps near the boundary. For completeness, the derivations of $\{X_{k+1}\}$ and $\{Y_{k+1}\}$ are outlined in Appendix~B.

\begin{lemma}\label{lem-2-2}
Let $\alpha \in (0, 2)$ and $r>0$ be the radius of the current ball. Assume that the ball jumps within the domain $\Omega$. Then the location of $(k+1)$-th ball, starting from the $k$-th ball, is given by
\begin{equation}\label{Jumptonext} 
X_{k+1}=X_k+J\cdot\begin{bmatrix}\cos\theta_1\\\sin\theta_1\cos\theta_2\\\cdots\cdots\\\sin\theta_1\cdots\sin\theta_{d-2}\sin\theta_{d-1}\end{bmatrix},
\end{equation}
where the jump distance $J$ is determined by
\begin{equation}
J(\xi,r,\alpha)=\frac{r}{\sqrt{I^{-1}(1-\xi;\alpha/2,1-\alpha/2)}},\;\;\;\xi\in (0,1),
\end{equation}
with $I^{-1}(\cdot\,;\,\cdot,\,\cdot)$ being the inverse  of the incomplete Beta function. Moreover, the random variable $Y_{k+1}$ associated with $X_k$ is given by
\begin{equation}\label{randvara} 
Y_{k+1}=X_k+\gamma \begin{bmatrix}\cos\theta_1\\
\sin\theta_1\cos\theta_2\\
\cdots\cdots\\
\sin\theta_1\cdots\sin\theta_{d-2}\sin\theta_{d-1}\end{bmatrix},
\end{equation}
where the random parameter $\gamma$ is given by
\begin{equation*}
\gamma :=\gamma(\xi,r,\alpha)= \xi^{1/\alpha} r,\;\;\;\xi\in (0,1).
\end{equation*}
\end{lemma}

With this, we can construct a Monte Carlo procedure based on the random sample as 
 \begin{equation}
 \label{Si}
S_{i}(\bx) = g(X_{m^{\ast}}^{i}) +  \sum_{k=0}^{m^{\ast}-1}\omega_{r_k^i}  f(Y_{k+1}^{i})W(X^i_k,Y^i_{k+1}),
 \end{equation}
 where the index $i$ represents the $i$-th trajectory and $X_0^i=\bx\in \Omega$.
 Accordingly, the solution to \eqref{fractional-laplace} admits the Monte Carlo approximation based on \eqref{th-2-1-1}
\begin{equation}\label{newalgo}
u(\bx)\approx \frac{1}{N}\sum_{i=1}^{N}S_{i}(\bx) =\frac{1}{N}\sum_{i=1}^{N}\Big[g(X^{i}_{m^{\ast}}) +  \sum_{k=0}^{m^{\ast}-1}\omega_{r_k^i}  f(Y^{i}_{k+1})W(X^i_k,Y^i_{k+1}) \Big],\;\;\;\bx\in \Omega.
 \end{equation}

 It is worth noting that the approximation \eqref{newalgo} is equivalent to the scheme developed in \cite{Sheng2023}. In particular, invoking \cite[Lemma 3.1]{Sheng2023} immediately yields that the estimator preserves unbiasedness, i.e., $\lim_{N\to\infty}N^{-1}\sum_{i=1}^NS_i(\bx)$ constitutes an unbiased estimator of $u(\bx)$. Denoting $\overline{S}(\bx)=N^{-1}\sum_{i=1}^NS_i(\bx)$, we arrive at the following error bound, whose proof is analogous to \cite[Theorem 3.1]{Sheng2023} and is therefore omitted for brevity.

\begin{cor}\label{thm-2.2}
For any $\varepsilon > 0$, assume that $u \in L_{\alpha}^1(\mathbb{R}^d)$, $f \in C^{\alpha+\varepsilon}(\Omega) \cap C(\overline{\Omega})$, $g \in L_{\alpha}^1(\Omega^c) \cap C(\Omega^c)$, and $\alpha \in (0, 2)$. Then, we have
\begin{equation}\label{errorbound}
\mathbb{E}\big[\overline{S}(\bx)-u(\bx)\big]^{2}
\le C (N^{-1}+\varepsilon^{2\alpha}).
\end{equation}
\end{cor}

From \eqref{errorbound}, we observe that we can
increase the number of samples  $N$  to reduce 
the error, which is costly due to the lower convergence rate.  To address this limitation, we propose in the forthcoming section a neural network based strategy to reduce the sampling cost.

\section{\bf{Two fractional neural walk-on-spheres methods}}\label{sec3}
{We first propose the FNWoS method, which enhances the simplified FWoS algorithm by integrating neural network surrogates to substantially reduce the number of trajectories required to reach a target accuracy. The key idea is to leverage a deep neural network (DNN) to approximate the output of the simplified FWoS estimator, thereby accelerating convergence and reducing the overall sampling cost.}
We further develop an enhanced method, coined as BFNWoS, by capping the maximum number of steps per sample path and incorporating a buffer strategy to eliminate the need to wait for the supervision signal to become highly accurate before training begins. We will show that the BFNWoS is more efficient than FNWoS for the case where $\alpha$ approaches $2$.

\subsection{{\bf Fractional neural walk-on-spheres method }}\label{sec31}
In this section, we present the detailed formulation of the FNWoS algorithm and outline its main components. We begin by recalling the Feynman–Kac representation for the fractional Poisson equation  \eqref{fractional-laplace}:
\begin{equation*} 
u(\bx) = \mathbb{E}_{X_{0}^{\alpha }=\bx}\big[g(X^{\alpha}_{\tau_{\Omega}})\big] 
        + \mathbb{E}_{X_{0}^{\alpha }=\bx}\Big[\int_{0}^{\tau_{\Omega}} f(X^{\alpha}_{s})\,{\rm d}s\Big],
\end{equation*}
where $\{X_{t}^{\alpha}\}_{t\geq0}$ 
where is a symmetric $\alpha$-stable L\'{e}vy process with initial point $X_{0}^{\alpha}=\bx \in \Omega$, and $\tau_{\Omega} = \inf \{ t>0 : X_{t}^{\alpha}\notin \Omega \}$.  

As described in the preceding section (see \eqref{newalgo}), the simplified FWoS scheme provides an efficient discretization of the probabilistic representation above. It yields a Monte Carlo estimator of the form:
\begin{equation}\label{FWoS}
u_{\text{FWoS}}(\bx_j)= \frac{1}{N}\sum_{i=1}^{N}\text{FWoS}^{i}(\bx_j),
 \end{equation}
 with individual trajectory contributions:
 \begin{equation*}\text{FWoS}^{i}(\bx_j)=g(X^{i,j}_{m^{\ast}}) +  \sum_{k=0}^{m^{\ast}-1}\omega_{r_k^{i,j}}  f(Y^{i,j}_{k+1})W(X^{i,j}_k,Y^{i,j}_{k+1}),
 \end{equation*}
 where $X_0^{i,j}=\bx_j\in \Omega$ is the initial point of the $i$-th trajectory, and the discrete paths $\{X_k^{i,j}\}_{k \ge 0}$ and $\{Y_k^{i,j}\}_{k \ge 1}$ are generated according to Lemma~\ref{lem-2-2}. The normalized weight function $\omega_{r_k^{i,j}}$ and the Kernel function $W$ are defined in Lemma \ref{lem2.1}. 
 
To accelerate the solution process, the FNWoS method introduces a neural nework surrogate model
$ v_{\theta}(\bx)$, parameterized by $\theta$. The  network takes spatial coordinates $\bx $ as input and produces an output vector whose dimension matches that of $ u(\bx) $. 
To construct the training dataset, we randomly sample $M$ input points $\{\bx_j\}_{j=1}^M$ and produce the corresponding reference outputs $u_{\mathrm{FWoS}}(\bx_j)$ using the simplified FWoS estimator \eqref{FWoS}.

The DNN surrogate $ v_{\theta}$ is then trained by minimizing the following empirical loss function:
\begin{equation} \label{lr-5-1-1}
\mathcal{L}_{\mathrm{FNWoS}}[v_\theta]
= \frac{1}{M} \sum_{j=1}^{M} \Big[ v_\theta(\bx_j)
    -u_{\text{FWoS}}(\bx_j) \Big]^2.
\end{equation}
This surrogate enables efficient evaluation of 
$u(\bx)$ at arbitrary query locations, effectively amortizing the cost of stochastic sampling across the domain and providing a smooth global approximation of the solution.

In practice, the simplified FWoS estimator already incorporates the boundary condition through the term 
$g$, so an explicit boundary loss is typically unnecessary. However, when the maximum step limit is small and trajectories may be truncated inside the domain, adding a boundary penalty can empirically improve stability and performance \cite{nam2024solving}. 

To this end, we partition the training data into interior points $\{\bx_j\}_{j=1}^{m_I}$ and boundary points $\{\bx_j\}_{j=1}^{m_B}$, and define a composite loss function: 
\begin{equation}\label{eq-loss-vr}
\hat{\mathcal{L}}_{\mathrm{FNWoS}}[v_\theta]
    := \frac{1}{m_I} \sum_{j=1}^{m_I} \Big[ v_\theta(\bx_j)
    -u_{\text{FWoS}}(\bx_j) \Big]^2+ \frac{\beta}{m_B} \sum_{j=1}^{m_B} \Big[ v_\theta(\bx_j)
    - g(\bx_j) \Big]^2,
\end{equation}
where $\beta> 0$ is a tunable penalty parameter that controls the relative weight of the boundary loss. This formulation ensures that for points lying on the boundary, the loss directly enforces consistency with the prescribed boundary condition. 
A summary of the complete FNWoS algorithm is presented in Algorithm \ref{fNWoS}.  

\begin{algorithm}[!ht]
\caption{Training of the FNWoS method}\label{fNWoS}
\small
\setlength{\algomargin}{0.1em}

\KwIn{number of iterations $T$, batch size $m_I$ and $m_B$ for interior and boundary points, 
source term $f$, boundary term $g$, order $\alpha$, stopping tolerance $\varepsilon$, 
number of trajectories $N$, boundary penalty parameter $\beta$}

\KwOut{optimized parameters $\theta$}

Initialize the neural network $v_\theta$;

    $\bx_\Omega = \{\bx_j\}_{j=1}^{m_I}\gets \text{sample $m_I$ points from } \Omega$  \tcc*{Sample points in $\Omega$}

    \For{$j \gets 1,\ldots,m_I$}{
    
        $\bx \gets \bx_j$ \tcc*{Current interior point}
        
        \For{$i \gets 1,\ldots,N$}{
        
            $S_i \gets 0$ ;
            $\bx^i \gets \bx$ ;
            $r \gets \text{dist}(\bx^i,\partial\Omega)$\tcc*{Compute distance to $\partial\Omega$}
            
            \While{$r > \varepsilon$}{
            
                $\gamma = r\,\zeta^{1/\alpha}, \quad \zeta \sim uniform(0,1)$;
                
                $y \gets \text{sample uniformly from }\partial B_{\gamma}(\bx^i)$;
                
                $S_i \gets S_i + \omega_r f(y) W(\bx^{i},y)$\tcc*{Estimate source}
                
                $J = r\left(I^{-1}(1-\omega;\alpha/2,1-\alpha/2)\right)^{-1/2},\ 
                \omega \sim uniform(0,1)$;
                
                $\Delta x \gets \text{sample uniformly from }\partial B_J(\bx^i)$;
                
                $\bx^i \gets \bx^i + \Delta x$ \tcc*{Walk to next point}
                
                $r \gets \text{dist}(\bx^i,\partial\Omega)$;
            }
            
            $\text{FWoS}^i(\bx_j) \gets S_i + g(\bx^i)$;
        }
        
        $u_{\text{FWoS}}(\bx_j) \gets \frac{1}{N}\sum_{i=1}^N \text{FWoS}^i(\bx_j)$;
    }
    
    $\by_\Omega \gets \{u_{\text{FWoS}}(\bx_j)\}_{j=1}^{m_I}$ \tcc*{Simplified FWoS estimates}
    
   \For{$k \gets 0, \ldots, T$}{
        $\bx_{\partial\Omega} \gets \text{sample $m_B$ points from $\partial\Omega$}$  \tcc*{Sample points in $\partial\Omega$}
    
        $\hat{\mathcal{L}}_{\text{I}} \gets \text{MSE}(v_\theta(\bx_\Omega), {\by}_{\Omega})$ \tcc*{Domain loss}
    
        $\hat{\mathcal{L}}_\text{B} \gets \text{MSE}(v_\theta(\bx_{\partial\Omega}), g(\bx_{\partial\Omega}))$ \tcc*{Boundary loss}
    
        $\hat{\mathcal{L}}_{\text{FNWoS}} \gets \hat{\mathcal{L}}_{\text{I}} + \beta \hat{\mathcal{L}}_\text{B}$ \tcc*{Combined loss}
    
        $\theta \gets \text{\texttt{step}}(\gamma, \nabla_{\!\theta} \hat{\mathcal{L}}_{\text{FNWoS}})$ \tcc*{SGD step}
}

\end{algorithm}

\begin{remark}{
We briefly explain why FNWoS outperforms simplified FWoS. From \eqref{lr-5-1-1}--\eqref{eq-loss-vr}, FNWoS determines $v_\theta$ by minimizing
$$
\mathcal{L}[v_\theta]
= \mathbb{E}\Big[(\,v_\theta(\bx)-u_\mathrm{FWoS}(\bx)\,)^2\Big],
$$
where we regard $u_\mathrm{FWoS}(\bx)$ as a zero-mean perturbation of the exact value $u(\bx)$.
Therefore, minimizing the above mean–squared loss implies that the DNN surrogate $v_\theta$ approximates the conditional expectation
\[
v_\theta\approx \mathbb E\big[u_\mathrm{FWoS}(\bx)|\bx\big]=u(\bx),
\]
i.e. it performs a regression to the mean of the scattered Monte Carlo samples produced by simplified FWoS.
Moreover, DNNs are known to possess a strong inductive bias toward low-frequency and smooth components of target functions {\rm\cite{xu2025overview}}. As a consequence, FNWoS projects the noisy pointwise estimates $u_\mathrm{FWoS}(\bx)$ onto a smooth function class, effectively filtering out stochastic fluctuations. This mechanism explains why FNWoS can achieve more accurate approximations.
}
\end{remark}
 
\subsection{{\bf {Buffered} fractional neural walk-on-spheres method}}

{Although the FNWoS algorithm described in Section 3.1 achieves better accuracy over the simplified FWoS under the same number of sampled trajectories, the maximum number of steps per trajectory increases markedly as the fractional order 
$\alpha$ becomes larger. To improve computational efficiency and scalability in such regimes, we further introduce the BFNWoS method, which incorporates a hybrid strategy that truncates each stochastic path after a limited number of steps and uses the current neural network surrogate to approximate the remaining contribution. In addition, a buffered supervision mechanism is employed to enable early stage training with partially accurate simplified FWoS outputs, without waiting for fully converged reference values. We now present the detailed formulation of the BFNWoS algorithm.}

\subsubsection{\bf{Stochastic path truncation}}
Although it is guaranteed that the discrete sequence of points $\{X_k\}_{k\ge 0}$ generated by the Markov process will eventually exit the given domain within a finite number of steps, it has been demonstrated in~\cite{Sheng2023} that the exit time increases substantially with both the dimensionality and the fractional order 
$\alpha$.  As a result, generating sufficient training data becomes increasingly costly in high-dimensional or large 
$\alpha$ settings. To mitigate this issue, we introduce a fixed maximum number of steps, denoted by 
$K$, for each stochastic trajectory. However, this truncation inevitably introduces a non-negligible approximation error when the trajectory remains inside the domain upon reaching step 
$K$.

 To address this issue, we decompose the simplified FWoS scheme \eqref{FWoS} into two components: the sum of the first 
$K$ terms, and the remaining terms which include the boundary contribution. Specifically, we write: 
\begin{equation}\label{FWoST}
    \mathrm{FWoS}^i(\bx_j)
    = \underbrace{\sum_{k=0}^{K-1} \omega_{r_k^{i,j}} f(Y_{k+1}^{i,j}) W(X_k^{i,j}, Y_{k+1}^{i,j})}_{\text{First $K$ terms}}
    + \underbrace{\sum_{k=K}^{m^*-1} \omega_{r_k^{i,j}} f(Y_{k+1}^{i,j}) W(X_k^{i,j}, Y_{k+1}^{i,j}) + g(X_{m^*}^{i,j})}_{\text{Remaining terms} \;=\; \mathrm{FWoS}^i(X_K^{i,j})},
\end{equation}
where $m^\ast= \inf\{\ell : X_\ell \notin \Omega\}$ is the true first exit step of the trajectory, and by construction $K\leq m^\ast$. Notably, the second part of this decomposition can be interpreted as the result of applying the simplified $\mathrm{FWoS}$ method starting from the intermediate point $X_K^{i,j}$.  
Motivated by this observation,
we propose to replace the remaining terms with a neural network surrogate to complete the truncated trajectory. To avoid ambiguity with the moving target, we distinguish the trainable network 
$v_{\theta}$ from a frozen network 
$v_{\bar{\theta}}$ used only for target generation. In practice, we initialize 
$\bar{\theta}\leftarrow \theta$ and then set $\bar{\theta}$
to the parameters from the  previous buffer refresh cycle, keeping 
$\bar{\theta}$ fixed when computing supervision signals.
\begin{figure}[!htbp]
\centering 
\includegraphics[width=0.70\textwidth]{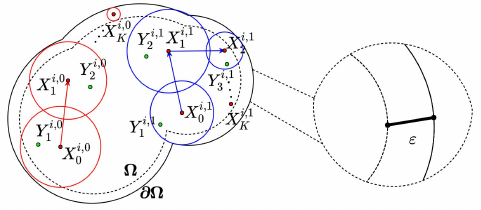} 
\caption{FNWoS with the maximum number of steps $K$. Red ball: center of the inscribed ball at the $K$-th step of the red trajectory lies in the $\varepsilon$-region; 
Blue ball: center of the inscribed ball at the $K$-th step of the blue trajectory lies in the interior of $\Omega$.
}
\label{NWoS}
\end{figure}

Accordingly, the simplified FWoS estimator in equation~\eqref{FWoS} can be reformulated under the truncation strategy described above. Specifically, we define the truncated FWoS estimator, denoted by $u_{\text{FWoS}_T}(\bx_j;\bar{\theta})$, as

\begin{equation}\label{eq-5-3-1}
u_{\text{FWoS}_T}(\bx_j;\bar{\theta})= \frac{1}{N}\sum_{i=1}^{N}\text{FWoS}_T^{i}(\bx_j;\bar{\theta}),
\end{equation}
with individual trajectory contributions given by:
\begin{equation}\label{eq-FWoST}
    \mathrm{FWoS}^i_{T}(\bx_j;\bar{\theta})
    = y^{\bx_j, v_{\bar{\theta}}} + \sum_{k=0}^{K-1} \omega_{r_k^{i,j}} f(Y_{k+1}^{i,j}) W(X_k^{i,j}, Y_{k+1}^{i,j}), 
\end{equation}
where $y^{\bx, v_{\bar{\theta}}}$ is defined as
\begin{equation} \label{eq:y_x_v}
    y^{\bx, v_{\bar{\theta}}} :=
    \begin{cases}
        v_{\bar{\theta}}(X_K^{i,j}), & d(X_K^{i,j}, \partial\Omega) > \varepsilon, \\
        g(X_K^{i,j}), & \text{otherwise},
    \end{cases}
\end{equation}
 to incorporate either the neural surrogate prediction or the exact boundary value, depending on whether the truncated point is sufficiently close to the boundary. Here, $d(X_K^{i,j}, \partial\Omega)$ denotes the distance from $X_K^{i,j}$ to the boundary $\partial\Omega$, and 
$\varepsilon$ is a small tolerance parameter. The two cases in \eqref{eq:y_x_v} are illustrated in Figure~\ref{NWoS}.

This leads naturally to a reformulated learning problem, in which the surrogate model 
$v_\theta$ is trained to match a truncated and network augmented estimator 
$u_{\mathrm{FWoS}_T}(\bx_j;\bar{\theta})$. Following the partition of samples into interior points 
$\{\bx_j\}_{j=1}^{m_I} \subset \Omega$ and boundary points 
$\{\bx_j\}_{j=1}^{m_B} \subset \partial\Omega$, the BFNWoS loss function is defined as
\begin{equation}\label{eq-Bloss-vr}
\hat{\mathcal{L}}_{\mathrm{BFNWoS}}[v_\theta]
    := \frac{1}{m_I} \sum_{j=1}^{m_I} \Big[ v_\theta(\bx_j)
    - u_{\mathrm{FWoS}_T}(\bx_j;\bar{\theta}) \Big]^2
    + \frac{\beta}{m_B} \sum_{j=1}^{m_B} \Big[ v_\theta(\bx_j)
    - g(\bx_j) \Big]^2,
\end{equation}
where $\beta > 0$ is a penalty parameter that balances the influence of interior and boundary terms.
\subsubsection{\bf{Data refinement based on neural caches}}
Inspired by \cite{li2023neural,nam2024solving}, we adopt buffered supervision strategy to alleviate the high computational cost associated with generating highly accurate training targets in the FNWoS method. Instead of waiting until all supervision signals reach high fidelity, we allow training to begin with approximate references, which are progressively refined throughout the training process.

To quickly initialize the buffer, we randomly sample $B$ input points $\{\bx_j\}_{j=1}^B\subset \Omega$ and compute the corresponding reference outputs 
$$\bY_j=\frac{1}{N_{init}}\sum_{i=1}^{N_{init}}\text{FWoS}_T^{i}(\bx_j;\bar{\theta})=\frac{1}{N_{init}}\sum_{i=1}^{N_{init}}(y^{\bx_j, v_{\bar{\theta}}} + \sum_{k=0}^{K_{init}-1} \omega_{r_k^{i,j}} f(Y_{k+1}^{i,j}) W(X_k^{i,j}, Y_{k+1}^{i,j}))$$ using \eqref{eq-5-3-1}  with a small number of trajectories $N_{init}$ and the current frozen $\bar\theta$. 
These pairs form the initial buffer $$\mathcal{B} = \{ ( \bx_j, \bY_j ) \}_{j=1}^{B}.$$
Here we employ a trajectory counter $fre[\bx_j]$ for each point $\bx_j$,  recording the total number of trajectories used to generate its current supervision signal $\bY_j$. During training, the buffer $\mathcal{B}$ is periodically updated every $L \in \mathbb{N}$ steps.  
At each update, we first randomly select a $\{(\bx^{update}_j, \bY^{update}_j)\}_{j=1}^{B_1}$ from $\mathcal{B}$ for refinement. For each $\bx^{update}_j$, we compute a new $\mathrm{FWoS}_T$ estimate $\by^{update}_j$ using an additional $N$ trajectories and update the outputs via  weighted averaging:
\[
\bY^{update}_j = \frac{fre[\bx^{update}_j]}{fre[\bx^{update}_j] + N} \bY^{update}_j + \frac{N}{fre[\bx^{update}_j] + N} \by^{update}_j,
\]
followed by incrementing the trajectory counter $$fre[\bx^{update}_j] = fre[\bx^{update}_j] + N.$$

Meanwhile, we randomly sample 
$\{(\bx^{\mathrm{replace}}_j, \bY^{\mathrm{replace}}_j)\}_{j=1}^{B_2}$ 
from $\mathcal{B}$ and replace them with newly sampled pairs 
$\{(\bx^{\mathrm{new}}_j, \by^{\mathrm{new}}_j)\}_{j=1}^{B_2}$, 
where each $\bx^{\mathrm{new}}_j$ is randomly sampled from $\Omega$ and 
$\by^{\mathrm{new}}_j$ is its $\mathrm{FWoS}_T$ estimate computed using $N$ trajectories.
This mechanism maintains a balance between diversity and accuracy in the training data.

The neural network parameters $\theta$ are optimized on the current buffer $\mathcal{B}$ by minimizing the loss in~\eqref{eq-Bloss-vr} using the Adam optimizer, and the buffer is refined and refreshed every $L$ training steps through the update–replace procedure described above. After each buffer refresh, we perform a hard update of the lagged network by setting 
$\bar{\theta}\leftarrow \theta$, and keep 
$\bar{\theta}$
 fixed during the subsequent target generation stage. For clarity, the BFNWoS method is summarized in Algorithm~\ref{al-5-1}. Let \( P_B \) denote the proportion of boundary samples in the training batch of size \( 2m \), and let \( P_{\mathcal{B}} \) denote the proportion of buffer samples refined in each update. 
The training batch thus consists of \( m_B = 2m P_B \) boundary samples and \( m_I = 2m (1 - P_B) \) interior samples. 
In each buffer update, \( B_1 = m P_{\mathcal{B}} \) samples are selected for refinement, \( B_2 = m (1 - P_{\mathcal{B}}) \) samples are replaced with newly generated trajectory pairs and buffer $\mathcal{B}$ of size $B=10m$. In particular, when the hyperparameter $L$ exceeds the total number of training iterations, BFNWoS reduces to FNWoS. Implementation details are provided in Appendix D.
\begin{algorithm}[!ht]
\caption{Training of the BFNWoS method}\label{al-5-1}
\small
\setlength{\algomargin}{0.1em}

\KwIn{$\text{truncated FWoS estimator}$ defined in \eqref{eq-5-3-1}, number of iterations $T$, warm-up iterations $C$, batch size $m_I$ and $m_B$ for interior and boundary points, batch size $B_1$ and $B_2$ for update points and replace points,
buffer $\mathcal{B}$ of size $B$, boundary function $g$, buffer update interval $L$, 
boundary penalty parameter $\beta$, number of trajectories $N$, number of init trajectories $N_{init}$, 
the maximum number of steps $K$ and $K_{init}$.}

\KwOut{optimized parameters $\theta$}

Initialize the networks $v_\theta$ ;

Initialize $\bar{\theta}\leftarrow \theta$ ;

$\bx_{\Omega} = \{\bx_j\}_{j=1}^{B} \gets \text{sample $B$ points from } \Omega$ \tcc*{Sample points in $\Omega$}

\For{$j \gets 1,\ldots,B$}{
    $\bY_j \gets \frac{1}{N_{init}}\sum_{i=1}^{N_{init}}{\text{FWoS}^i_{T}}(\bx_j; \bar{\theta})$;
    
    $fre[\bx_j] = N_{init}$ \tcc*{Truncated FWoS estimator with $K_{init}$ and $N_{init}$}
}

$\mathcal{B} \gets \text{initialize with } (\bx_{\Omega}, \bY_\Omega=\{\bY_j\}_{j=1}^B)$\tcc*{Initial buffer with $K_{init}$ and small $N_{init}$}

\For{$k \gets 0, \ldots, T$}{
    \If{$k \bmod L = 0$ and $k > C$}{
        
        $\bar{\theta}\leftarrow \theta$ ;
        
        $\bx^{update}_{\mathcal{B}_1} = \{\bx_j\}_{j=1}^{B_1} \gets \text{sample $B_1$ points from } \mathcal{B}$ ;

        \For{\text{each} $\bx_j \in \bx^{update}_{\mathcal{B}_1}$}{
            $\by_j^{update} \gets \frac{1}{N}\sum_{i=1}^N {\text{FWoS}^i_{T}}(\bx_j; \bar{\theta})$\tcc*{Truncated FWoS estimator with $K$ and $N$}
            
            $\bY_j
            \gets
            \frac{fre(\bx_j)}{fre(\bx_j)+N}\,\bY_j
            +
            \frac{N}{fre(\bx_j)+N}\,\by_j^{update}$;
            
            $fre(\bx_j) \gets fre(\bx_j) + N$\tcc*{Update estimates}
        }

        $\bx^{replace}_{\mathcal{B}_2} \gets \text{sample $B_2$ points from } \mathcal{B}$ ;
        
        $\bx^{new}_{\Omega} = \{x_j\}_{j=1}^{B_2} \gets \text{sample $B_2$ points from } \Omega$ ;
        
        \For{\text{each} $\bx_j \in \bx^{new}_{\Omega}$}{
            $\bY_j^{new} \gets \frac{1}{N}\sum_{i=1}^N {\text{FWoS}^i_{T}}(\bx_j; \bar{\theta})$\tcc*{Truncated FWoS estimator with $K$ and $N$}
            
            $fre(\bx_j) \gets N$ \tcc*{Replace}
        }

        $\mathcal{B} \gets \text{randomly replace}~
        (\bx^{replace}_{\mathcal{B}_2}, \bY^{replace}_{\mathcal{B}_2})
        ~\text{with}~
        (\bx^{new}_{\Omega}, \bY^{new}_{\Omega}=\{\bY_j^{new}\}_{j=1}^{B_2})$ ;
        }    
    
    $\bx_{\partial\Omega} \gets \text{sample $m_B$ points from } \partial\Omega$ \tcc*{Sample points in $\partial\Omega$}
    
    $(\bx_\mathcal{B}, \bY_\mathcal{B}) \gets \text{sample $m_I$ pairs from } \mathcal{B}$ \tcc*{Sample points in $\mathcal{B}$}
    
    $\hat{\mathcal{L}}_{\text{I}} \gets \text{MSE}(v_\theta(\bx_\mathcal{B}), \bY_\mathcal{B})$ \tcc*{Domain loss}
    
    $\hat{\mathcal{L}}_\text{B} \gets \text{MSE}(v_\theta(\bx_{\partial\Omega}), g(\bx_{\partial\Omega}))$ \tcc*{Boundary loss}
    
    $\hat{\mathcal{L}}_{\text{BFNWoS}} \gets \hat{\mathcal{L}}_{\text{I}} + \beta \hat{\mathcal{L}}_\text{B}$ \tcc*{Combined loss}
    
    $\theta \gets \text{\texttt{step}}(\theta, \nabla_\theta \hat{\mathcal{L}}_{\text{BFNWoS}})$ \tcc*{SGD step}
}
\end{algorithm}

  \section{Numerical Results}\label{sec4}
In this section, we present representative numerical results that demonstrate the efficiency of the proposed methods, with an emphasis on Algorithms~\ref{fNWoS} and \ref{al-5-1} for solving~\eqref{fractional-laplace}. When reporting the numerical results, FWoS is short for the simplified FWoS.
In all subsequent examples, we define the relative $\ell^2$-error as
$$\text{Relative $\ell^2$-error} = \big( \| u_{\ast} - u \|_{\ell^2(\Omega)}\big)/\big(\| u \|_{\ell^2(\Omega)}\big),\;\;\;\| u \|_{\ell^2(\Omega)}:=\Big(\frac{1}{n} \sum_{j=1}^nu^2(\bx_j)\Big)^{\frac12},$$
where $u_\ast(\bx)$ denotes the numerical solution computed by MC-fPINNs, FWoS, FNWoS, or BFNWoS, and $u(\bx)$ is the exact solution. The error is evaluated over $n=10^5$ uniformly sampled points $\{\bx_j\}_{j=1}^n$ in $\Omega$. Detailed parameter settings can be found in Appendix~D.

\begin{example}\label{ex4.1} 
{\rm({\bf Benchmark problem).}}
We first consider the problem \eqref{fractional-laplace} on the unit ball $\Omega = \mathbb{B}^d_1=\{\bx\in\mathbb{R}^d:|\bx| < 1\}$
subject to homogeneous boundary conditions. The exact solution is 
$u(\bx) = (1 - |\bx|^2)_+^{1 + \frac{\alpha}{2}}$, where $a_+ := \max\{a,0\}$, and the right-hand side source function is taken as 
\begin{equation*}
\begin{split}
f(\bx) = 2^\alpha \Gamma\Big(\frac{\alpha}{2} + 2\Big)   \Gamma\Big(\frac{\alpha + d}{2}\Big)   \Gamma\Big(\frac{d}{2}\Big)^{-1} 
        \Big(1 - \big(1 + \frac{\alpha}{d}\big)|\bx|^2\Big).
\end{split}
\end{equation*}

\end{example}

\begin{figure}[htbp]
  \centering
  \begin{subfigure}[b]{0.33\textwidth}
    \includegraphics[width=\linewidth]{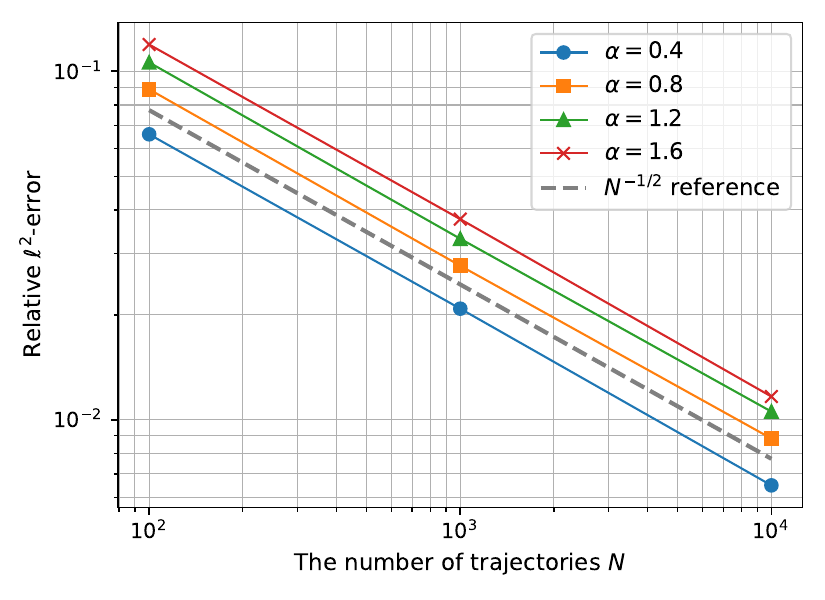}
    \caption{Convergence rate of FWoS}
    \label{fig:poisson_nd_wos}
  \end{subfigure}\hspace{-4pt}
  \hfill 
  \begin{subfigure}[b]{0.33\textwidth}
    \includegraphics[width=\linewidth]{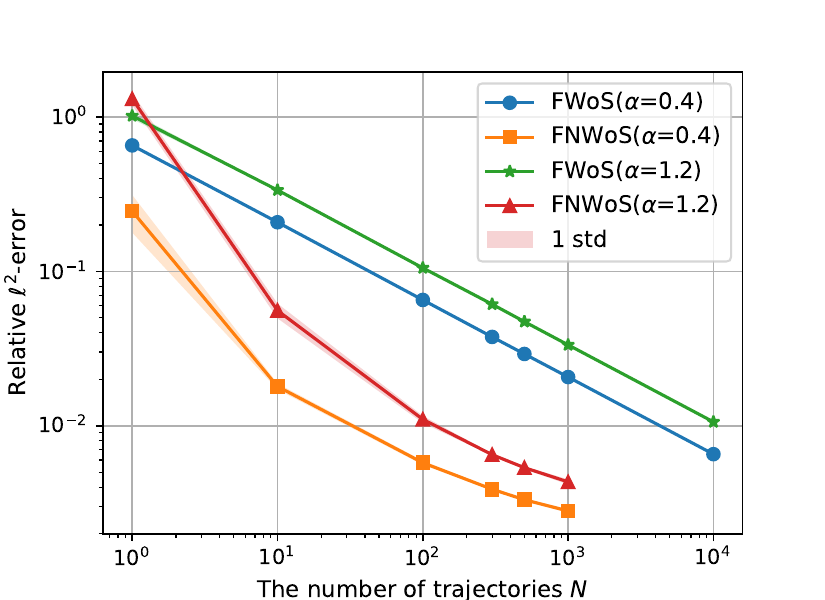}
    \caption{ FWoS vs FNWoS}
    \label{fig:WoS_NWoS_ball_compare}
  \end{subfigure}\hspace{-4pt}
  \hfill 
  \begin{subfigure}[b]{0.33\textwidth}
    \includegraphics[width=\linewidth]{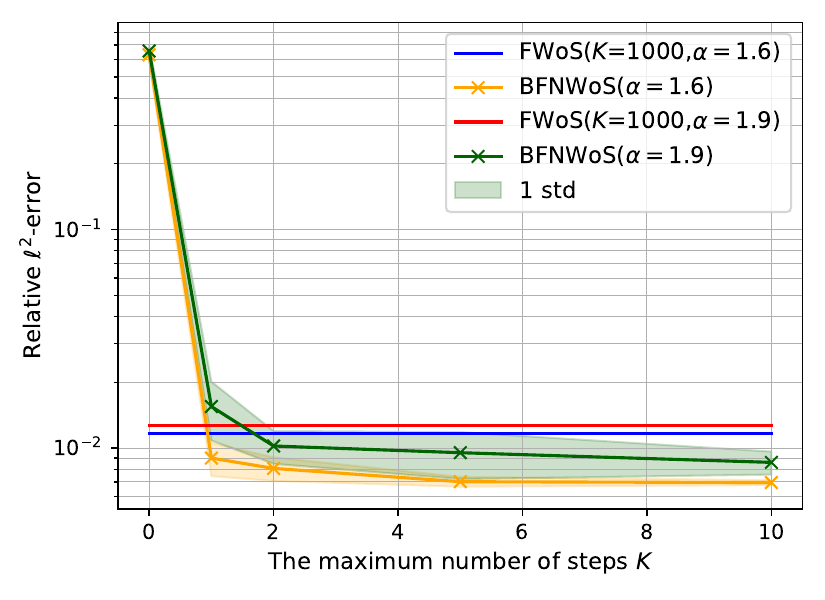}
    \caption{FWoS vs BFNWoS}
    \label{fig:max_step}
  \end{subfigure}
  \caption{The relative $\ell^2$-errors of FWoS, FNWoS and BFNWoS for Example 4.1 on the 10D unit ball. (a): The relative $\ell^2$-errors of FWoS versus the parameter $N$ with various $\alpha$. 
  (b): Comparison of the relative $\ell^2$-errors of FWoS and FNWoS with respect to trajectories $N$ for $\alpha = 0.4$ and $1.2$.  
(c): Comparison of the relative $\ell^2$-errors of FWoS and BFNWoS with respect to the maximum number of steps $K$. The error of FWoS is shown as a horizontal line at $\alpha = 1.6$ and $1.9$. The colored lines and shaded regions correspond to mean values and one-standard-deviation bands of relative $\ell^2$-errors.}
  \label{fig:WoS_NWoS_ball_compare and fig:max_step}
\end{figure}

We first employ the FWoS method to solve  \eqref{fractional-laplace} on the 10-dimensional unit ball to verify the effectiveness of our algorithm, where the maximum number of steps is fixed at $K=1000$. Figure~\ref{fig:WoS_NWoS_ball_compare and fig:max_step}(a) shows the relative $\ell^2$-errors for fractional orders $\alpha \in \{0.4, 0.8, 1.2, 1.6\}$. It can be observed that, for all the fractional orders considered, the convergence rate is $\mathcal{O}(N^{-1/2})$, which is consistent with the theoretical result in Corollary~\ref{thm-2.2}. We next compare the FWoS method with FNWoS in Figure~\ref{fig:WoS_NWoS_ball_compare and fig:max_step}(b), with emphasis on the dependence of accuracy on the trajectory number $N$. We perform five independent runs for each fractional order. As shown in Figure~\ref{fig:WoS_NWoS_ball_compare and fig:max_step}(b), for all tested $\alpha$, the relative $\ell^2$-error decreases as $N$ increases. More importantly, FNWoS $(K=1000,\,N=100,\,\beta=10,\,\text{iterations}=40000)$ attains higher accuracy than FWoS $(K=1000,\,N=10000)$ while using far fewer trajectories. 

Since Markov trajectories exit the domain after sufficiently many steps, the required step count grows sharply with increasing dimension and as $\alpha\to 2$, which reduces training efficiency. We therefore compare BFNWoS with a step cap against FWoS at $\alpha=1.6$ and $1.9$. We run five independent trials per fractional order. Figure~\ref{fig:WoS_NWoS_ball_compare and fig:max_step}(c) reports the one standard deviation band of the relative $\ell^2$-error, where BFNWoS uses $K_{init}=1000$, $N_{{init}}=1$, $N=100$, $\beta=10$, $L=100$, $C=1$, $\text{iterations}=40000$, whereas FWoS uses $N=10000$. Once the maximum step count exceeds a threshold, the relative $\ell^2$-error plateaus. In the classical FWoS (see also Fig.~9 of \cite{Sheng2023}), the average step count exceeds $400$ for large $\alpha$, whereas our BFNWoS attains comparable accuracy with   $1$ or $2$ steps.

\begin{table}[H]
\centering
\caption{The relative $\ell^2$-errors of MC-fPINNs, FWoS, FNWoS and BFNWoS for Example 4.1 on the 50D unit ball with various $\alpha$. Under the provided hyperparameter setup, the solutions obtained by FNWoS with $N=100$ and BFNWoS with $N=100$ are much more accurate than that of FWoS with $N=400$.}
\label{tab:wos_comparison_swos2}
\small 
\setlength{\tabcolsep}{5pt} 
\begin{tabularx}{0.99\textwidth}{
  >{\RaggedRight}X
  *{6}{>{\centering\arraybackslash}X}
}
\toprule
\textbf{Method}
& \textbf{$\alpha=0.4$}
& \textbf{$\alpha=0.8$}
& \textbf{$\alpha=1.2$}
& \textbf{$\alpha=1.9$} \\
\midrule
MC-fPINNs   
& 1.15e-1 & 2.30e-1 & 5.36e-2  & \textbf{4.01e-2} \\
FWoS        
& 3.87e-2 & 5.02e-2  & 5.76e-2 & 6.67e-2 \\
FNWoS
& \textbf{3.79e-2} & 4.41e-2 & 4.13e-2 &  4.64e-2 \\
BFNWoS
& 3.87e-2 & \textbf{4.00e-2} &  \textbf{3.87e-2} & 4.62e-2 \\
\bottomrule
\end{tabularx}
\end{table}

    \begin{figure}[!ht]
    \centering
    \begin{subfigure}[b]{0.48\textwidth}
        \includegraphics[width=\textwidth]{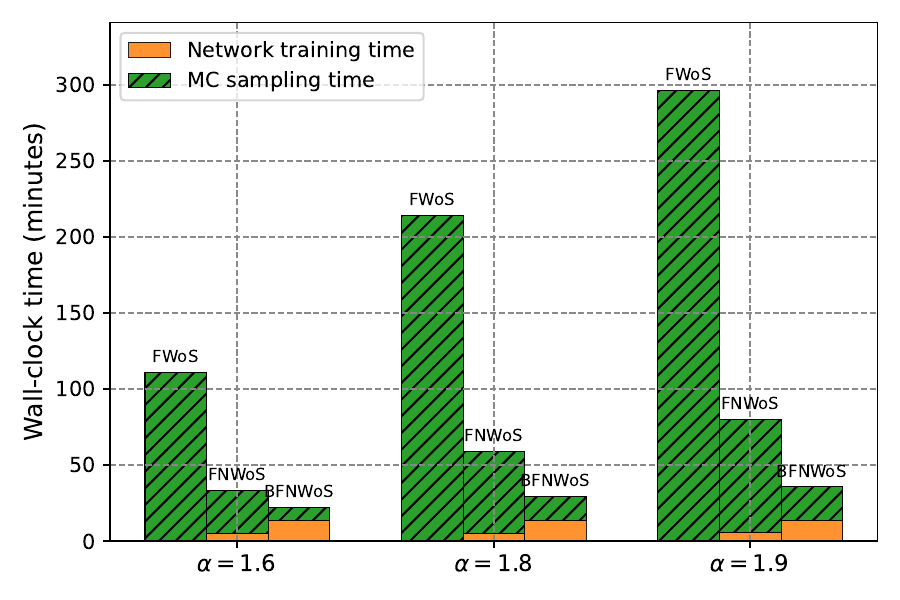 }
        \caption{Computational time breakdown}
    \end{subfigure}\hspace{6pt}
    \begin{subfigure}[b]{0.48\textwidth}
        \includegraphics[width=\textwidth]{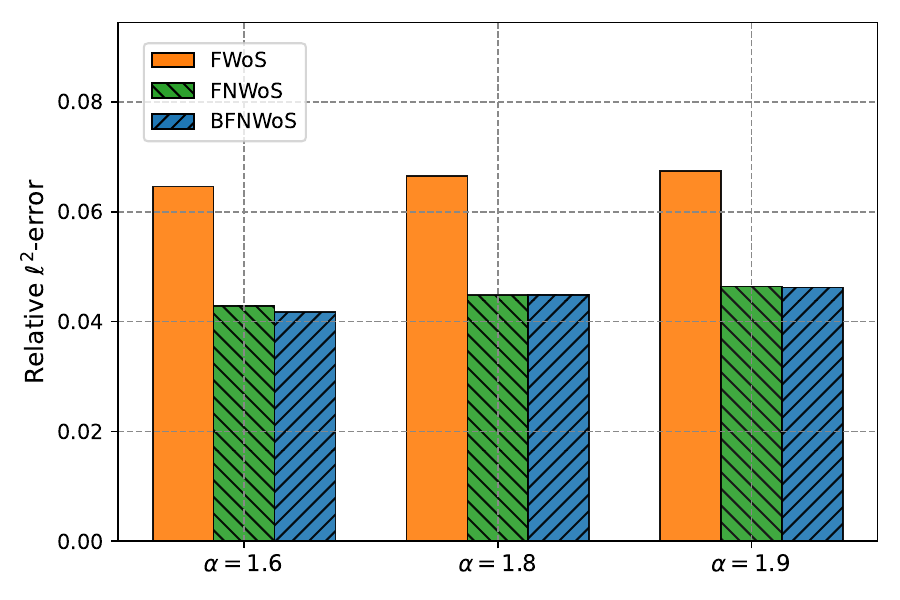}
        \caption{Prediction accuracy}
    \end{subfigure}
    \caption{
Computational time breakdown and corresponding prediction accuracy of FWoS, FNWoS, and BFNWoS for Example 4.1 on the 50D unit ball. 
(a) Wall-clock time averaged over five independent runs, decomposed into Monte Carlo sampling time used to evaluate the FWoS estimator or to generate supervision targets and network training time used to refine supervision targets and train the network. FWoS has no training stage, so its training time is zero. 
(b) Corresponding prediction accuracy measured by relative $\ell^2$-error, obtained from the same runs as in (a).
}
    \label{fig:time and rl2}
\end{figure}

{We next compare the relative $\ell^2$-errors of MC-fPINNs, FWoS, FNWoS, and BFNWoS. As reported in Table \ref{tab:wos_comparison_swos2}, FNWoS and BFNWoS achieve noticeably higher accuracy than FWoS using fewer trajectories, although their errors remain slightly larger than those of MC-fPINNs on this benchmark in some cases. We note that this test case involves a regular domain and a solution with relatively high regularity. These conditions are generally favorable for achieving high accuracy with MC-fPINNs. Nevertheless, the methods proposed in this paper are directly applicable to irregular domains and to general boundary conditions, including nonhomogeneous ones. We will demonstrate this broader applicability and improved robustness in the following numerical examples.}

We also report the computational time and corresponding prediction accuracy of the last three methods in Figure~\ref{fig:time and rl2}(a) and Figure~\ref{fig:time and rl2}(b), respectively, averaged over five independent runs.
The results show that the proposed FNWoS and BFNWoS methods become more efficient than the FWoS method as $\alpha$ increases. 
In particular, the required step count grows markedly with larger 
$\alpha$, so FNWoS spends substantially more time to generate the training set. In contrast, the initial training set is constructed more quickly by BFNWoS than by FNWoS. The benefit of its truncation strategy becomes more pronounced for larger $\alpha$. Indeed, BFNWoS can complete the training phase before FNWoS finishes data generation. For smaller $\alpha$, where the required step count is relative small, this advantage diminishes, and a nontrivial share of runtime is spent on buffer updates and replacements, which reduces the overall efficiency. This performance is consistent with our expectations.

\begin{example}\label{ex4.1.2} 
{\rm({\bf Low-regularity solution on the unit disk with homogeneous BCs).}}
We next consider the unit disk $\Omega=\mathbb{B}_1^2$ with the discontinuous source term $f(r,\theta)=\mathbf{1}_{\{|\theta|<\pi/2\}}$, where $(r,\theta)$ are polar coordinates. Then, the exact solution is given in the form reported in~\cite{ainsworth2017aspects}:
\begin{equation}
    \begin{aligned}
        u(r,\theta) = &2^{-\alpha} \left( \frac{(1 - r^2)_+^{\alpha/2}}{\Gamma(1 + \alpha/2)^2} \right) \biggl\{ \frac{1}{2}\\
        &+ \sum_{\substack{n \geq 0 \\ \ell \geq 1 \\ \ell \text{ odd}}} \frac{(-1)^{(\ell+1)/2 + n + 1} (2n + \alpha/2 + \ell + 1) \cos(\ell \theta) r^\ell P_n^{(\alpha/2, \ell)}(2r^2 - 1)}{\pi (n + \ell/2)(\alpha/2 + 1)\binom{n + \alpha/2 + \ell/2 + 1}{n + \ell/2}\binom{\alpha/2 + n}{n}} \biggr\}, 
    \end{aligned}
\end{equation}
where $P_{n}^{(a,b)}(\cdot)$ denotes the Jacobi polynomial of degree $n$ with parameters $(a,b)$, and
$\binom{a}{b}$ denotes the generalized binomial coefficient defined by
$\binom{a}{b}=\frac{\Gamma(a+1)}{\Gamma(b+1)\Gamma(a-b+1)}$, with $\Gamma(\cdot)$ being the Gamma function.

\end{example}

\begin{figure}[!ht]
    \centering
    \begin{subfigure}[b]{0.425\textwidth}
        \includegraphics[width=\textwidth]{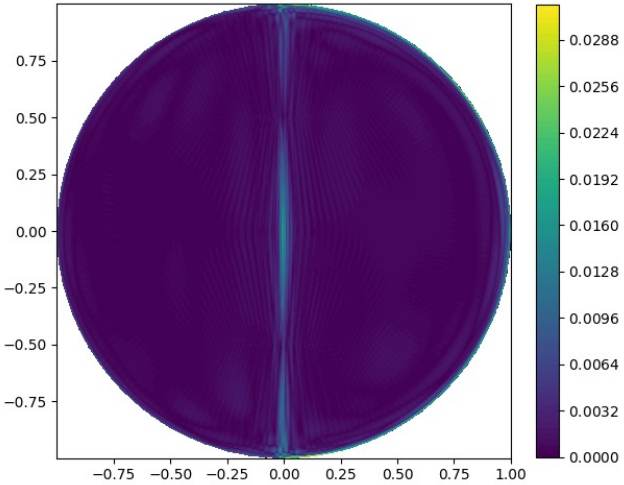}
        \caption{$\alpha=0.5$}
    \end{subfigure}\hspace{6pt}
    \begin{subfigure}[b]{0.41\textwidth}
        \includegraphics[width=\textwidth]{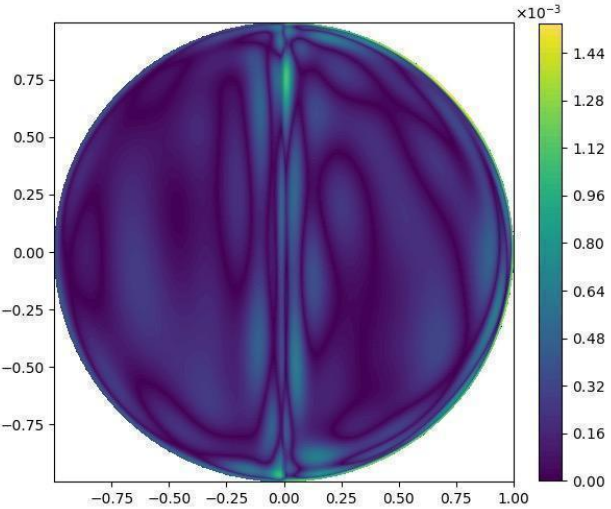}
        \caption{$\alpha=1.5$}
    \end{subfigure}
    \caption{The accuracy of FNWoS for Example~4.2 on the 2D unit disk with a low-regularity solution. 
(a): absolute error for $\alpha = 0.5$. (b): absolute error for $\alpha = 1.5$.}
    \label{fig:2d_ball_contour}
\end{figure}

{The solution of this problem exhibits low regularity on the unit disk. To demonstrate the performance of our methods, we report the relative $\ell^2$-errors in Table~\ref{tab:pde_2d_ball} and pointwise errors of the FNWoS method for $\alpha = 0.5$ and $1.5$ in Figure~\ref{fig:2d_ball_contour}. We observe that FNWoS and BFNWoS yield more accurate numerical results than FWoS, even with fewer trajectories. Notably, MC-fPINNs is evaluated using the same parameter settings as in~\cite{guo2022monte},
with detailed configurations provided in Appendix~D.
In comparison, FNWoS and BFNWoS consistently yield more accurate solutions. 
This highlights the advantage of our methods in low-regularity regimes, where Monte Carlo approximation errors become a bottleneck for MC-fPINNs.}

\begin{table}[H]
\centering
\caption{The relative $\ell^2$-errors of MC-fPINNs, FWoS, FNWoS and BFNWoS for Example 4.2 on the 2D unit disk with various $\alpha$. Under the provided hyperparameter setup, the solutions obtained by FNWoS with $N=1000$ and BFNWoS with $N=100$ are 
much more accurate than that of FWoS with $N=10000$.}
\label{tab:pde_2d_ball}
\footnotesize
\begin{tabularx}{0.98\textwidth}{ 
  >{\RaggedRight}X 
  >{\centering}X   
  *{4}{>{\centering\arraybackslash}X} 
}
\toprule
\textbf{Method} 
& \multicolumn{1}{c}{\textbf{$\alpha=0.5$}} 
& \multicolumn{1}{c}{\textbf{$\alpha=1.5$}} \\
\midrule
 MC-fPINNs   & 2.15e-2 &4.96e-3     \\
FWoS         & 7.89e-3  & 1.15e-2 \\
FNWoS  & \textbf{5.99e-3} & \textbf{1.45e-3} \\
BFNWoS  & 7.65e-3 & 2.67e-3 \\
\bottomrule
\end{tabularx}
\end{table}

\begin{example}\label{ex4.2} {\rm({\bf High-dimensional problem with nonhomogeneous BCs).}}
We next consider problem \eqref{fractional-laplace} on $\Omega=(0,1)^d$ with a source term
\begin{equation}
    f(\bx) = \frac{d\, 2^\alpha \Gamma\!\left(\frac{\alpha+d}{2}\right) \Gamma\!\left(\frac{\alpha+3}{2}\right)}
                 {\Gamma\!\left(\frac{3}{2}\right)\Gamma\!\left(\frac{d}{2}\right)}
                 (1 + |\bx|^2)^{-\frac{\alpha+3}{2}}
                 {}_2F_1\!\bigg(\frac{\alpha+3}{2}, -\frac{\alpha}{2}; \frac{d}{2}; \frac{|\bx|^2}{1+|\bx|^2}\bigg),
\end{equation}
and global nonhomogeneous boundary condition $g(\bx)=d(1+|\bx|^2)^{-3/2}$. 
This problem admits the exact solution $u(\bx)=d(1+|\bx|^2)^{-3/2}$ {\rm(cf.~\cite{sheng2020fast})}.
\end{example}

\begin{figure}[!ht]
    \centering
    \begin{subfigure}[b]{0.42\textwidth}
        \includegraphics[width=\textwidth]{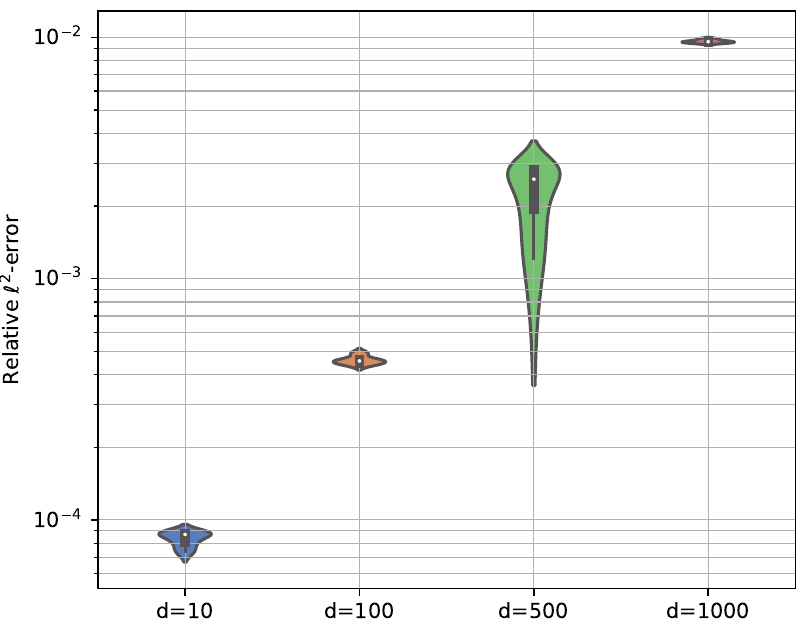}
        \caption{$\alpha=0.4$}
    \end{subfigure}\hspace{6pt}
    \begin{subfigure}[b]{0.42\textwidth}
        \includegraphics[width=\textwidth]{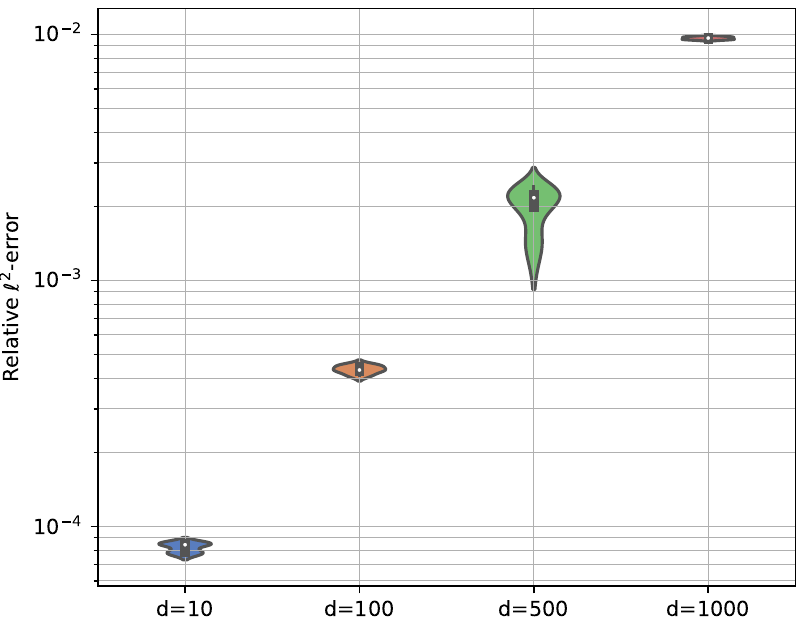}
        \caption{$\alpha=0.8$}
    \end{subfigure}  
    \caption{The relative $\ell^2$-errors of the FNWoS for Example 4.3 across dimensions $d = 10, 100, 500, 1000$. Each violin plot represents the distribution of errors over five independent experiments.
    }
    \label{fig:high dimension compare}
\end{figure}

Table~\ref{tab:pde_comparison_ann} reports the relative $\ell^2$-errors of the three proposed methods on a $d$-dimensional hypercube domains $\Omega=(0,1)^d$ with $d=10$. We observe that FNWoS and BFNWoS achieve higher accuracy than FWoS, which is consistent with the previous observations. {Here the MC-fPINNs is evaluated using the same parameter settings as in~\cite{guo2022monte},
with detailed configurations provided in Appendix~D.
In comparison, FNWoS and BFNWoS consistently yield more accurate solutions. This highlights the advantage of our methods in problems with nonhomogeneous boundary conditions, where the structural assumptions used in MC-fPINNs \cite{guo2022monte} no longer apply, and Monte Carlo approximation errors become a bottleneck for MC-fPINNs.}
We further test the performance of FNWoS method on higher-dimensional problems with fractional orders $\alpha=0.4$ and $0.8$. To account for stochastic variability, we run five independent trials per fractional order. Figure~\ref{fig:high dimension compare} plots the distribution of the relative $\ell^2$-error for $\alpha=0.4$ and $0.8$ at $d = 10, 100, 500$ and 1000, which indicates that FNWoS is effective in very high dimensions.

\begin{table}[H]
\centering
\caption{The relative $\ell^2$-errors of MC-fPINNs, FWoS, FNWoS and BFNWoS for Example 4.3 on $\Omega=(0,1)^{10}$ with various $\alpha$. Under the provided hyperparameter setup, the solutions obtained by FNWoS with $N=300$ and BFNWoS with $N=100$ are much more accurate than that of FWoS with $N=10000$.}
\label{tab:pde_comparison_ann}
\small 
\setlength{\tabcolsep}{5pt} 
\begin{tabularx}{0.99\textwidth}{
  >{\RaggedRight}X
  *{6}{>{\centering\arraybackslash}X}
}
\toprule
\textbf{Method}
& \textbf{$\alpha=0.4$}
& \textbf{$\alpha=0.8$}
& \textbf{$\alpha=1.2$}
& \textbf{$\alpha=1.9$} \\
\midrule
MC-fPINNs   
& 1.17e-3 & 1.38e-3 & 1.79e-3 & 2.91e-3 \\
FWoS        
& 5.41e-3 & 3.84e-3 & 2.52e-3 & 1.19e-3 \\
FNWoS  
& 8.58e-5 & 8.31e-5 & 9.56e-5 & 8.46e-5 \\
BFNWoS 
& \textbf{8.37e-5} & \textbf{8.25e-5} & \textbf{8.26e-5} & \textbf{8.17e-5} \\
\bottomrule
\end{tabularx}
\end{table}

\begin{example}\label{ex4.3}  
{\rm({\bf Source problem on irregular domains).}} 
Finally, we consider an irregular domain $\Omega\subset\mathbb{R}^d$ with a  source term {\rm(cf.~\cite{sheng2020fast})}
\begin{equation}
    f(\bx)=\frac{ 2^\alpha\Gamma((\alpha+d)/2)}{\Gamma(d/2)}~_1F_1\Big(\frac{\alpha+d}{2};{d/2};-|\bx|^2\Big),
\end{equation}
and the boundary condition $g(\bx)=e^{-|\bx|^2}$. This problem admits the exact solution $u(\bx) = e^{-|\bx|^2}$.

\end{example}

\begin{figure}[!ht]
    \centering
    \begin{subfigure}[b]{0.45\textwidth}
        \includegraphics[width=\textwidth]{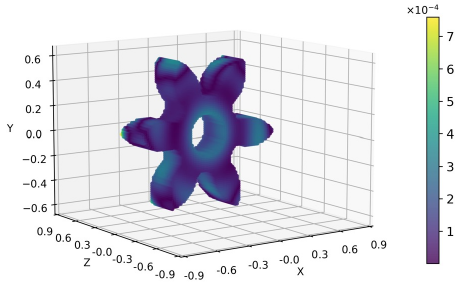}
        \caption{$\alpha=0.4$}
    \end{subfigure}
    \begin{subfigure}[b]{0.45\textwidth}
        \includegraphics[width=\textwidth]{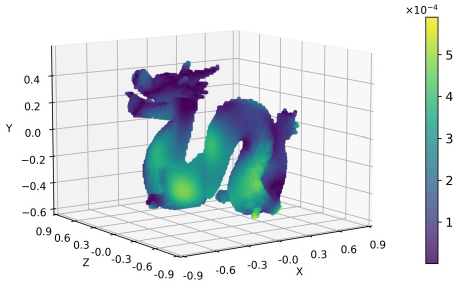}
        \caption{$\alpha=0.8$}
    \end{subfigure}
        \begin{subfigure}[b]{0.45\textwidth}
        \includegraphics[width=\textwidth]{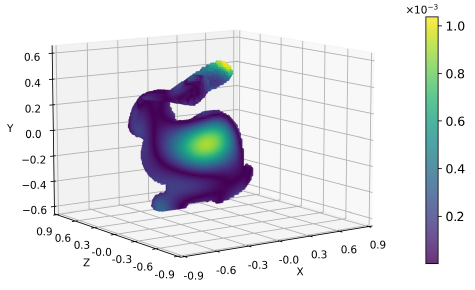}
        \caption{$\alpha=1.2$}
    \end{subfigure} 
        \begin{subfigure}[b]{0.45\textwidth}
        \includegraphics[width=\textwidth]{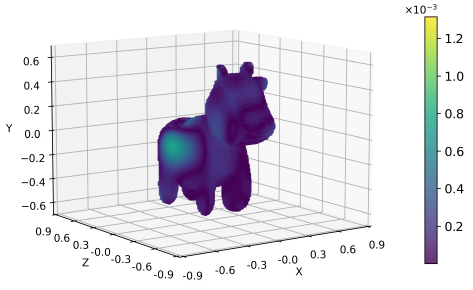}
        \caption{$\alpha=1.6$}
    \end{subfigure} 
    \caption{The accuracy of FNWoS and BFNWoS  for Example~4.4 on 3D animal-shaped domains with the fabricated solutions $u(\bx)=e^{-|\bx|^2}$. (a)--(c): absolute error obtained by FNWoS; (d): absolute error obtained by BFNWoS.}
    \label{fig:dragon_domain}
\end{figure}

We solve \eqref{fractional-laplace} on animal shaped domains and plot pointwise errors for the different domains in Figure~\ref{fig:dragon_domain}. Each domain $\Omega$ is represented by a signed distance function (SDF) realized in this experiment by a neural field based on the instant-NGP~\cite{li2023neural,muller2022instant}. We test the FNWoS method on irregular domains for $\alpha\in\{0.4,0.8,1.2\}$ (see Figures~\ref{fig:dragon_domain}(a)–(c)), and BFNWoS method for $\alpha=1.6$ (see Figure~\ref{fig:dragon_domain}(d)). We pre-generate a pool of boundary points at initialization, since boundary sampling is inefficient on such animal-shaped domains. At each iteration, we randomly select a subset for training to improve computational efficiency.

\section{Conclusion}\label{sec5}
We propose a simplified FWoS scheme and developed its neural-accelerated extensions, FNWoS and BFNWoS, for solving high-dimensional fractional Poisson equations on irregular domains. By simplifying the weight and sampling components in classical FWoS (cf. \cite{Sheng2023}), amortizing Monte Carlo sampling through a global neural surrogate, and further improving efficiency via truncated trajectories and buffered supervision, the proposed methods achieve favorable accuracy cost trade-offs and strong scalability in challenging regimes. Extensive numerical experiments demonstrate the robustness and efficiency of our approach across a wide range of dimensions and geometries. As an important direction for future work, we will extend the proposed framework to time-dependent problems and develop corresponding walk-on-spheres solvers for fractional parabolic equations.

\section*{Acknowledgements}

 The first author is supported by the NSF of China (under grant numbers 92270115, 12071301), the Shanghai Municipal Science and Technology Commission (No. 20JC1412500), and Henan Academy of Sciences. Hao Wu is supported by the NSF of China (under grant number 12571463). Fanhai Zeng is supported by the NSF of China (under grant number 12171283, 12531020).

\bigskip

 \renewcommand{\theequation}{A.\arabic{equation}}
\noindent{\bf \large Appendix A: Proof of Lemma \ref{lem2.1} } \label{AppendixA}

\medskip
Recall the integral form for the solution of \eqref{fractional-laplace} on $\B^d_r$(cf. \cite[Lemma 2.1]{Sheng2023}):
\begin{equation}\label{eq-lem-Tr00}\begin{aligned}
u(\boldsymbol{x})&=\int_{\mathbb{R}^d\setminus\mathbb{B}_r^d}g(\boldsymbol{z})
P_r(\boldsymbol{x},\boldsymbol{z})\dx[\boldsymbol{z}]
+\int_{\mathbb{B}_r^d}f(\boldsymbol{y})G_r(\boldsymbol{x},\boldsymbol{y})\dx[\boldsymbol{y}].
\end{aligned}\end{equation}

According to \cite[Lemma A.5]{MR3461641}, the function $P_r(\bx,\bz)$ serves as the probability density function for $g(\bz)$, satisfying
$\dint_{\mathbb{R}^{d}\setminus \B^d_{r}} P_{r}(\bx,\bz)\,\d \bz=1$, 
which implies that the first term in \eqref{eq-lem-Tr00} can be equivalently expressed as a conditional expectation:
\begin{equation}\label{eq-lem-2-2-3}
    \int_{\mathbb{R}^d \setminus \mathbb{B}_r^d} g(\boldsymbol{z}) P_r(\boldsymbol{x}, \boldsymbol{z})\dx[\boldsymbol{z}]=\mathbb{E}_{P_r}[g(Z)].
\end{equation}

Next, we reformulate the second part of \eqref{eq-lem-Tr00} into an expectation form. 
Direct computation shows that the Green’s function in \eqref{Gr} can be 
written as
\begin{equation*}
    \begin{aligned}
G_r(\boldsymbol{x}, \boldsymbol{y})&=\widetilde{C}_d^\alpha |\boldsymbol{y} - \boldsymbol{x}|^{\alpha -d} \int_0^{\varrho(\boldsymbol{x}, \boldsymbol{y})} \frac{t^{\frac{\alpha}{2} - 1}}{(t + 1)^{\frac{d}{2}}} \, \mathrm{d}t
\qquad (t=1/m-1)\\
    &=\widetilde{C}_d^\alpha|\boldsymbol{y}-\boldsymbol{x}|^{\alpha-d}\int_{\varrho^*(\boldsymbol{x},\boldsymbol{y})}^1\left(\frac{1}{m}-1\right)^{\alpha/2-1}\left(\frac{1}{m}\right)^{-\frac{d}{2}}\frac{1}{m^2}\dx[m] \\
        &=\widetilde{C}_d^\alpha|\boldsymbol{y}-\boldsymbol{x}|^{\alpha-d}\int_{\varrho^*(\boldsymbol{x},\boldsymbol{y})}^1m^{\frac{d-\alpha}{2}-1}(1-m)^{\alpha/2-1}\dx[m]\\
        &=\underbrace{\widetilde{C}_d^\alpha B{\Big(\frac{d-\alpha}2,\frac\alpha2\Big)}|\boldsymbol{y}-\boldsymbol{x}|^{\alpha-d}}_{\Q_r(\boldsymbol{x}, \boldsymbol{y})}
        \underbrace{\left[1-I{\Big(\varrho^*(\boldsymbol{x},\boldsymbol{y});\frac{d-\alpha}2,\frac\alpha2\Big)}\right]}_{W(\boldsymbol{x}, \boldsymbol{y})}\\
        &=W(\boldsymbol{x}, \boldsymbol{y})\Q_r(\boldsymbol{x}, \boldsymbol{y})
        =\omega_r W(\boldsymbol{x}, \boldsymbol{y})\frac{\Q_r(\boldsymbol{x}, \boldsymbol{y})}{\omega_r},
    \end{aligned}
\end{equation*}
where $I(\cdot\,;\,\cdot,\cdot)$ denotes the incomplete Beta function and
$\omega_r=\int_{\mathbb{B}_r^d}  \Q_r(\boldsymbol{x}, \boldsymbol{y}) \, \mathrm{d}\boldsymbol{y}
=\widetilde{{C}}_d^\alpha B(\frac{d-\alpha}{2},\frac{\alpha}{2})\int_{\mathbb{B}_r^d}|\boldsymbol{y} - \boldsymbol{x}|^{\alpha - d} \d\by=\frac{r^\alpha B(\frac{d-\alpha}{2},\frac{\alpha}{2})}{\alpha 2^{\alpha-1} \Gamma^2(\alpha/2)}$.
Let
\begin{equation*}
 Q_{r}(\bx,\by)= \frac{\Q_{r}(\bx,\by)}{\omega_r}
= \frac{\Gamma(d/2)r^\alpha}{2\alpha\pi^{d/2}} |\boldsymbol{y} - \boldsymbol{x}|^{\alpha -d}.
\end{equation*}
Then, we have $\dint_{\B^d_{r}} Q_{r}(\bx,\bz)\,\d \bz=1$ and
\begin{equation}\label{rho}\begin{aligned}
\int_{\mathbb{B}_r^d}f(\boldsymbol{y})G_r(\boldsymbol{x},\boldsymbol{y})\dx[\boldsymbol{y}]
=\omega_r\int_{\mathbb{B}_r^d}f(\boldsymbol{y})W(\boldsymbol{x},\boldsymbol{y})Q_r(\boldsymbol{x},\boldsymbol{y})\dx[\boldsymbol{y}]
=\omega_r\mathbb{E}_{Q_r}\big[f(Y)W(\boldsymbol{x},Y)\big].
\end{aligned}\end{equation}
From \eqref{eq-lem-Tr00}, \eqref{eq-lem-2-2-3}, and \eqref{rho}, we derive 
\eqref{eq-lem-2-2-1}.
This ends the proof.

\bigskip

\noindent{\bf \large Appendix B: Proof of Lemma~\ref{lem-2-2}} \label{AppendixB}
 \renewcommand{\theequation}{B.\arabic{equation}}
 
 \medskip

We first generate the sample sequence $\{X_k\}$, where the key step is to determine the jump distances, which can be obtained via the inverse transform sampling method. For any $r>0$ and  
$\boldsymbol{x}\in \mathbb{B}_r^d$, the following formula holds
\begin{equation*}
    \int_{\mathbb{R}^n\setminus\mathbb{B}_r^d}P_r(\boldsymbol{x},\boldsymbol{z})\,\d\boldsymbol{z}=1.
\end{equation*}
Direct calculations show
$|\boldsymbol{x}-\boldsymbol{z}|^2=|\boldsymbol{x}|^2+|\boldsymbol{z}|^2-2|\boldsymbol{x}||\boldsymbol{z}|\cos\theta=\rho^2+|\boldsymbol{x}|^2-2|\boldsymbol{x}|\rho\cos\theta$ and   consequently
\begin{equation*}
	\begin{aligned}
		&\int_{r<|\boldsymbol{z}|<\gamma}P_r(\boldsymbol{x},\boldsymbol{z}){\d\boldsymbol{z}}\\
        &=\widehat{C}_d^\alpha(r^2-|\bx|^2)^{\alpha/2}2\pi\prod_{k=1}^{ d-3}\int_0^\pi\sin^k\theta \,\d\theta\int_r^\gamma\int_0^\pi\frac{\rho^{ d-1}\sin^{ d-2}\theta \;\d\theta \d\rho}{(\rho^2-r^2)^{\alpha/2}(\rho^2+|\bx|^2-2\rho|\bx|\cos\theta)^{ d/2}},\\
	\end{aligned}
\end{equation*}
where we used $\bar{r} = r/|\bx| > 1$, $\bar{\rho} = \rho/|\bx| > 1$. According to \cite[(A.25)]{MR3461641} and \cite[(A.29)]{MR3461641}, the following identity holds:
\begin{equation*}\begin{split}
  \pi \prod_{k=1}^{d-2} \int_0^\pi \sin^k \theta \,~ \mathrm{d}\theta = \frac{\pi^{d/2}}{\Gamma(d/2)},
\end{split}\end{equation*}
and 
\begin{equation*}\begin{split}
  \int_0^\pi \frac{\sin^{d-2} \theta}{(\bar{\rho}^2 + 1 - 2\bar{\rho} \cos \theta)^{d/2}} \, \d\theta= \frac{1}{\bar{\rho}^{d-2}(\bar{\rho}^2 - 1)} \int_0^\pi \sin^{d-2} \theta \, \d\theta.
\end{split}\end{equation*}
Therefore, one obtains
\begin{equation*}
		\int_{r<|\boldsymbol{z}|<\gamma}P_r(\boldsymbol{x},\boldsymbol{z}){\d\boldsymbol{z}}\\
		=\frac{2\sin(\pi\alpha/2)}{\pi}(\overline{r}^2-1)^{\alpha/2}\int_{\overline{r}}^{\overline{\gamma}} \frac{\overline{\rho} }{(\overline{\rho}^2-\overline{r}^2)^{\alpha/2}(\overline{\rho}^2-1)}\d\overline{\rho}.
\end{equation*}
Next, by performing the substitution $\left(t - \frac{t}{\bar{r}^2}\right) = \left(\frac{\bar{\rho}}{\bar{r}}\right)^2 - 1$, we arrive at
\begin{equation*}
	\begin{aligned}
		&\int_{r<|\boldsymbol{z}|<\gamma}P_r(\boldsymbol{x},\boldsymbol{z}){\d\boldsymbol{z}}=\frac{\sin(\pi\alpha/2)}{\pi}\int_{0}^{\frac{\overline{\gamma}^2-\overline{r}^2}{\overline{r}^2 - 1}} \frac{1}{(1+t)t^{\alpha/2}} \, \d t.
	\end{aligned}
\end{equation*}
Letting $\widehat{t} = \frac{1}{1+t}$, one further obtains
\begin{equation*}
	\begin{aligned}
		\int_{r<|\boldsymbol{z}|<\gamma}P_r(\boldsymbol{x},\boldsymbol{z}){\d\boldsymbol{z}}&=\frac{\sin(\pi\alpha/2)}{\pi}\int_{\frac{\overline{r}^2-1}{\overline{\gamma}^2-1}}^1\widehat{t}^{-\alpha/2}(1-\widehat{t})^{\alpha/2-1}\:\mathrm{d}\widehat{t},\\
        &=\frac{\sin(\pi\alpha/2)}{\pi}B(\alpha/2,1-\alpha/2)\Big(1-I\Big(\frac{\overline{r}^2-1}{\overline{\gamma}^2-1},\alpha/2,1-\alpha/2\Big)\Big).
	\end{aligned}
\end{equation*}
Finally, noting the identity $\sin(\alpha\pi/2)B(\alpha/2,1-\alpha/2)=\pi$, we conclude 
\begin{equation*}
	\begin{aligned}
		\int_{r<|\boldsymbol{z}|<\gamma}P_r(\boldsymbol{x},\boldsymbol{z}){\d\boldsymbol{z}}=1-I\Big(\frac{\overline{r}^2-1}{\overline{\gamma}^2-1},\alpha/2,1-\alpha/2\Big):=\omega\in (0,1),
	\end{aligned}
\end{equation*}
where $\frac{\overline{r}^2-1}{\overline{\gamma}^2-1}=\frac{{r}^2-|\bx|^2}{{\gamma}^2-|\bx|^2}$. 
In the special case where $\bx$ coincides with both the origin of the coordinate system and the center of the ball $\B_{r}^d$, this reduces to
 \begin{equation*}
	\begin{aligned}
		\int_{r<|\boldsymbol{z}|<\gamma}P_r(\boldsymbol{x},\boldsymbol{z})\d\boldsymbol{z}=1-I\Big(\frac{{r}^2}{{\gamma}^2},\alpha/2,1-\alpha/2\Big):=\omega\in (0,1).
	\end{aligned}
\end{equation*}

We now turn to the sampling of $\{Y_k\}$. It is known that the relation $\int_{0<|\boldsymbol{y}|<r}Q_r(\boldsymbol{x},\boldsymbol{y})\dx[\boldsymbol{y}]=1$.  For all $ \gamma\in(0,r)$ by applying the spherical coordinate transformation, we can obtain that
\begin{equation*}
\begin{aligned}
 \int_{0<|\boldsymbol{y}|<\gamma}Q_r(\boldsymbol{x},\boldsymbol{y})
 \dx[\boldsymbol{y}]
 &=\frac{\alpha\Gamma(d/2)}{2\pi^{d/2}r^\alpha}\int_{0<|\boldsymbol{y}|<\gamma}
 {|\boldsymbol{y}|}^{\alpha-d}\dx[\boldsymbol{y}]
 =\Big(\frac{\gamma}{r}\Big)^{\alpha}:=\xi\in(0,1),
\end{aligned}
\end{equation*}
which indicates that the sample distance of $\boldsymbol{y}$ can be obtained as $\gamma = \xi^{1/\alpha} r$.
This completes the proof.

\bigskip

\noindent{\bf \large Appendix C: Proof of Theorem~\ref{th-2-1}} \label{AppendixC}
 \renewcommand{\theequation}{C.\arabic{equation}}
 
 \medskip 
 
We partition the trajectory of the stochastic process $X_{t}^{\alpha}$ within the domain $\Omega$ by segmenting it based on the escape times from each inscribed small ball. More precisely, we define the sequence $0 =\tau_{0}<\tau_{1}<\tau_{2}<\cdots<\tau_{m^{\ast}} = \tau_{\Omega}$, where $\tau_{k}$ $(0<k<m^{\ast})$ denotes the first exit time from the $k$-th small ball $\B_{r_k}^d:=\B_{r_k}^d(\bx_k)$ with radius $r_k$ and centered at $\bx_k$, while $\tau_{\Omega}$ represents the first exit time from the domain $\Omega$. Since the fractional Laplacian $(-\Delta)^{\frac{\alpha}2}$ and its associated Green’s function are translation invariant, we may, without loss of generality, shift the coordinate system so that the current center $\bx_k$ of $\B_{r_k}^d(\bx_k)$ is moved to the origin. Consequently, based on \eqref{Feynmankac}, the solution can be reformulated as
\begin{equation*}
\begin{split}
u(\bx) &= \mathbb{E}_{X_0 = \boldsymbol{x}} \left[ g(X_{\tau_\Omega}^\alpha) \right] + \mathbb{E}_{X_0 = \boldsymbol{x}} \left[ \int_0^{\tau_\Omega} f(X_s^\alpha) \, \d s \right] \\
&= \mathbb{E}_{X_0 = \boldsymbol{x}} \left[ g(X_{\tau_\Omega}^\alpha) \right] + \mathbb{E}_{X_0 = \boldsymbol{x}} \left[ \int_{\tau_0}^{\tau_1} f(X_s^\alpha) \, \d s + \int_{\tau_1}^{\tau_2} f(X_s^\alpha) \, \d s + \dots + \int_{\tau_{m^*-1}}^{\tau_{m^*}} f(X_s^\alpha) \, \d s \right] \\
&= \mathbb{E}_{X_0 = \boldsymbol{x}} \left[ g(X_{\tau_{m^*}}^\alpha) \right] + \mathbb{E}_{X_0 = \boldsymbol{x}} \left[ \sum_{k=0}^{m^*-1} \int_{\tau_k}^{\tau_{k+1}} f(X_s^\alpha) \, \d s \middle| X_{\tau_k}^\alpha \right].
\end{split}
\end{equation*}
By \eqref{Feynmankac} and \eqref{eq-lem-2-2-1} when $f(x)=0$, we obtain that
\begin{equation*}
\mathbb{E}_{X_0^\alpha=\boldsymbol{x}}\left[g(X_{\tau_m*}^\alpha)\right]=\int_{\mathbb{R}^d\setminus\mathbb{B}_{r_{m^*-1}}^d}g(\boldsymbol{z})P_r(\boldsymbol{x},\boldsymbol{z})\mathrm{d}\boldsymbol{z}.
\end{equation*}
Next, extracting the $k$-th term from the summation, by \eqref{Feynmankac} and \eqref{eq-lem-2-2-1} when $g(\bx)=0$, we have
\begin{equation*}\begin{aligned}
\mathbb{E}_{X_0 = \boldsymbol{x}} \Big[ \int_{\tau_k}^{\tau_{k+1}} f(X_s^\alpha) \, ds \big| X_{\tau_k}^\alpha = \boldsymbol{x}_k \Big] = \omega_{r_k} \int_{\mathbb{B}_{r_k}^d} f(\boldsymbol{y})W(\bx_k, \boldsymbol{y})Q_{r_k}(\boldsymbol{x}_k,\boldsymbol{y}) \, \d \by \\
= \mathbb{E}_{Q_{r_k}} \Big[ \omega_{r_k}f(Y_{k+1}) W(\boldsymbol{x}_k, Y_{k+1}) \big| X_{\tau_k}^\alpha = \boldsymbol{x}_k \Big], \quad Y_{k+1} \in \mathbb{B}_{r_k}^d.
\end{aligned}\end{equation*}
Finally, combining the above results, we can obtain the desired result.

\bigskip

\noindent{\bf \large Appendix D: Implementation Details} \label{implementation details}
 \renewcommand{\theequation}{B.\arabic{equation}}
 
 \medskip

Some implementation details, ablation study and parameter settings are follows.
\begin{itemize}
\item The proposed algorithms were implemented using the PyTorch framework, with their details presented in pseudocode in Algorithms \ref{fNWoS} and \ref{al-5-1}. The experimental evaluations were carried out on an A800 GPU platform. We choose a feedforward neural network with residual connections. 
Model training employed the Adam optimizer with an initial learning rate of $0.001$. At each iteration, the learning rate was multiplicatively decayed by a factor of $0.01^{1/T}$, where $T$ is the total number of iterations.
\item In practical computation, we choose appropriate values of $L$ and the maximum number of steps $K$ to match the specific learning problem and the method. In practice, when we set $L$ larger than the total number of iterations $T$ and $K_{init}$ is large, BFNWoS reduces to FNWoS, and the two algorithms are unified in one framework.
\item For MC-fPINNs, training was performed using the Adam optimizer with an initial learning rate of $0.001$, which was multiplied by $0.1$ at iterations $25000$, $30000$, $35000$, and $40000$ via a MultiStepLR scheduler. For Examples~4.1 and~4.2, the MC-fPINNs enforces the homogeneous boundary condition by embedding it directly into the neural network architecture, i.e., the trial solution is constructed as $(1 - |\bx|^2)_+ \cdot u_{\mathrm{NN}}(\bx)$, where $u_{\mathrm{NN}}$ denotes the output of the neural network.
\item Table~\ref{tab:ablation} shows the ablation study for the replacement strategy in BFNWoS. Without replacement, the error is much larger and increases quickly as $\alpha$ grows. With replacement, the error stays small for all $\alpha$ values. This shows that the replacement strategy is important for good performance.

\begin{table}[H]
\centering
\caption{Ablation study of replacement in BFNWoS for Example~4.1 on the 10D unit ball with various $\alpha$.}
\label{tab:ablation}
\footnotesize
\setlength{\tabcolsep}{9pt}
\begin{tabular}{@{} l c *{6}{c} @{}}
\toprule
\textbf{Method}  & $\alpha=0.4$ & $\alpha=0.8$ & $\alpha=1.2$ & $\alpha=1.6$ & $\alpha=1.8$ & $\alpha=1.9$ \\
\midrule
BFNWoS w/o replacement    & 1.58e-2 & 5.14e-2 & 9.37e-2 & 1.56e-1 & 1.83e-1 & 1.99e-1 \\
BFNWoS w/ replacement  & \textbf{3.75e-3} & \textbf{5.47e-3} & \textbf{6.49e-3} & \textbf{7.11e-3} & \textbf{8.71e-3} & \textbf{1.05e-2} \\
\bottomrule
\end{tabular}
\end{table}

\begin{table}[H]
\centering
\caption{The parameter settings for the results in Table~\ref{tab:wos_comparison_swos2}.}
\footnotesize
\begin{tabularx}{0.95\textwidth}{
  >{\RaggedRight}X
  *{4}{>{\centering\arraybackslash}X}
}
\toprule
& \textbf{FWoS} 
& \textbf{FNWoS} 
& \textbf{BFNWoS} \\
\cmidrule(lr){2-2} \cmidrule(lr){3-3} \cmidrule(lr){4-4}
$\alpha$
&\{0.4,0.8,1.2,1.9\} &\{0.4,0.8,1.2,1.9\} &\{0.4,0.8,1.2,1.9\} \\
width, depth &  & 256, 6 & 256, 6 \\
Activation &  & GELU & GELU \\
Learning rate &  & $1\times10^{-3}$ & $1\times10^{-3}$ \\
$m,P_B,P_\mathcal{B}$ &  & $2^{13},0.1,0.6$ & $2^{13},0.1,0.6$ \\
$\varepsilon$ & $1\times10^{-4}$ & $1\times10^{-4}$ & $1\times10^{-4}$  \\
$K_{init},N_{init}$  &  & 10000,100 & 10000,30 \\
$K,N$ & 10000,400 & 10000,100 & 1,100 \\
$\beta$ &  & 10 & 10 \\
$L,C$ &  & 1000000,20000 & 100,20000 \\
Iterations &  & 40000 & 40000 \\
\bottomrule
\end{tabularx}
\end{table}

\begin{table}[H]
\centering
\caption{The parameter settings for the results in Figure~\ref{fig:time and rl2}.}
\footnotesize
\begin{tabularx}{0.95\textwidth}{
  >{\RaggedRight}X
  *{4}{>{\centering\arraybackslash}X}
}
\toprule
& \textbf{FWoS} 
& \textbf{FNWoS} 
& \textbf{BFNWoS} \\
\cmidrule(lr){2-2} \cmidrule(lr){3-3} \cmidrule(lr){4-4}
$\alpha$
&\{1.6,1.8,1.9\} &\{1.6,1.8,1.9\} & \{1.6,1.8,1.9\} \\
width, depth &  & 256, 6 & 256, 6 \\
Activation &  & GELU & GELU \\
Learning rate &  & $1\times10^{-3}$ & $1\times10^{-3}$ \\
$m,P_B,P_\mathcal{B}$ &  & $2^{13},0.1,0.6$ & $2^{13},0.1,0.6$ \\
$\varepsilon$ & $1\times10^{-4}$ & $1\times10^{-4}$ & $1\times10^{-4}$  \\
$K_{init},N_{init}$  &  & 10000,100 & 10000,30 \\
$K,N$ & 10000,400 & 10000,100 & 1,100 \\
$\beta$ &  & 10 & 10 \\
$L,C$ &  & 1000000,1 & 100,20000 \\
Iterations &  & 40000 & 40000 \\
\bottomrule
\end{tabularx}
\end{table}

\begin{table}[H]
\centering
\caption{The parameter settings for the results in Table~\ref{tab:pde_2d_ball}.}
\footnotesize
\begin{tabularx}{0.95\textwidth}{
  >{\RaggedRight}X
  *{4}{>{\centering\arraybackslash}X}
}
\toprule
& \textbf{FWoS} 
& \textbf{FNWoS} 
& \textbf{BFNWoS} \\
\cmidrule(lr){2-2} \cmidrule(lr){3-3} \cmidrule(lr){4-4}
$\alpha$
&\{0.5,1.5\} &\{0.5,1.5\} & \{0.5,1.5\} \\
width, depth &  & 256, 6 & 256, 6 \\
Activation &  & GELU & GELU \\
Learning rate &  & $1\times10^{-3}$ & $1\times10^{-3}$ \\
$m,P_B,P_\mathcal{B}$ &  & $2^{13},0.1,0.6$ & $2^{13},0.1,0.6$ \\
$\varepsilon$ & $1\times10^{-20}$ & $1\times10^{-20}$ & $1\times10^{-20}$  \\
$K_{init},N_{init}$  &  & 1000,1000 & 1000,1 \\
$K,N$ & 1000,10000 & 1000,1000 & 1,100 \\
$\beta$ &  & 1 & 1 \\
$L,C$ &  & 1000000,1 & 100,1 \\
Iterations &  & 150000 & 150000 \\
\bottomrule
\end{tabularx}
\end{table}

\begin{table}[H]
\centering
\caption{The parameter settings for the results in Table~\ref{tab:pde_comparison_ann}.}
\footnotesize
\begin{tabularx}{0.95\textwidth}{
  >{\RaggedRight}X
  *{4}{>{\centering\arraybackslash}X}
}
\toprule
& \textbf{FWoS} 
& \textbf{FNWoS} 
& \multicolumn{2}{c}{\textbf{BFNWoS}} \\
\cmidrule(lr){2-2} \cmidrule(lr){3-3} \cmidrule(lr){4-5}
$\alpha$
&\{0.4,0.8,1.2,1.9\} &\{0.4,0.8,1.2,1.9\} & \{0.4,0.8,1.2\} & \{1.9\} \\
width, depth &  & 256, 6 & 256, 6 & 256, 6 \\
Activation &  & GELU & GELU & GELU \\
Learning rate &  & $1\times10^{-3}$ & $1\times10^{-3}$ & $1\times10^{-3}$ \\
$m,P_B,P_\mathcal{B}$ &  & $2^{13},0.1,0.6$ & $2^{13},0.1,0.6$ & $2^{13},0.1,0.6$ \\
$\varepsilon$ & $1\times10^{-4}$ & $1\times10^{-4}$ & $1\times10^{-4}$ & $1\times10^{-4}$ \\
$K_{init},N_{init}$  &  & 1000,300 & 1000,1 & 1000,10 \\
$K,N$ & 1000,10000 & 1000,300 & 1,100 & 1,100 \\
$\beta$ &  & 5000 & 5000 & 5000 \\
$L,C$ &  & 1000000,1 & 100,1 & 100,1 \\
Iterations &  & 100000 & 100000 & 100000 \\
\bottomrule
\end{tabularx}
\end{table}

\begin{table}[H]
\centering
\caption{The parameter settings for the results in Figure~\ref{fig:high dimension compare}.}
\footnotesize
\begin{tabularx}{0.95\textwidth}{
  >{\RaggedRight}X
  *{4}{>{\centering\arraybackslash}X}
}
\toprule
& \multicolumn{4}{c}{\textbf{FNWoS}}\\
\cmidrule(lr){2-2} \cmidrule(lr){3-3} \cmidrule(lr){4-4} \cmidrule(lr){5-5}
$d$ &10 &100 &500 &1000\\
$\alpha$
 &\{0.4,0.8\} & \{0.4,0.8\} & \{0.4,0.8\} & \{0.4,0.8\} \\
width, depth   & 256, 6 & 512, 6 & 512, 6& 512, 6 \\
Activation   & GELU & GELU & GELU & GELU \\
Learning rate   & $1\times10^{-3}$ & $1\times10^{-3}$ & $1\times10^{-3}$ & $1\times10^{-3}$ \\
$m,P_B,P_\mathcal{B}$   & $2^{13},0.1,0.6$ & $2^{13},0.1,0.6$ & $2^{14},0.1,0.6$ & $2^{14},0.1,0.6$ \\
$\varepsilon$ & $1\times10^{-4}$ & $1\times10^{-4}$ & $1\times10^{-4}$ & $1\times10^{-4}$ \\
$K_{init},N_{init}$ & 1000,1000 & 1000,1000 & 1000,1 & 1000,1 \\
$K,N$ & 1000,1000 & 1000,1000 & 1000,1 & 1000,1 \\
$\beta$  & 5000 & 5000 & 5000 & 5000 \\
$L,C$  & 1000000,1 & 1000000,1 & 1000000,1 & 1000000,1 \\
Iterations & 100000 & 100000 & 100000 & 100000 \\
\bottomrule
\end{tabularx}
\end{table}

\begin{table}[H]
\centering
\caption{The parameter settings for the results in Figure~\ref{fig:dragon_domain}.}
\footnotesize
\begin{tabularx}{0.95\textwidth}{
  >{\RaggedRight}X
  *{2}{>{\centering\arraybackslash}X}
}
\toprule
& \textbf{FNWoS} 
& \textbf{BFNWoS} \\
\cmidrule(lr){2-2} \cmidrule(lr){3-3}
$\alpha$
 &\{0.4,0.8,1.2\} & \{1.6\} \\
width, depth   & 256, 6 & 256, 6 \\
Activation   & GELU & GELU \\
Learning rate   & $1\times10^{-3}$ & $1\times10^{-3}$ \\
$m,P_B,P_\mathcal{B}$   & $2^{12},0.001,0.6$ & $2^{12},0.001,0.6$ \\
$\varepsilon$ & $1\times10^{-4}$ & $1\times10^{-4}$ \\
$K_{init},N_{init}$ & 1000,1000 & 1000,1\\
$K,N$ & 1000,1000 & 1,100\\
$\beta$  & 500 & 500 \\
$L,C$  & 1000000,1 & 100,1 \\
Iterations & 40000 & 40000\\
\bottomrule
\end{tabularx}
\end{table}

\begin{table}[H]
\centering
\caption{The parameter settings of MC-fPINNs.}
\footnotesize
\begin{tabularx}{0.95\textwidth}{
  >{\RaggedRight}X
  *{4}{>{\centering\arraybackslash}X}
}
\toprule
& \textbf{Example 4.1} 
& \textbf{Example 4.2} 
& \textbf{Example 4.3} \\
\cmidrule(lr){2-2} \cmidrule(lr){3-3} \cmidrule(lr){4-4}
$\alpha$
&\{0.4,0.8,1.2,1.9\} &\{0.5,1.5\} & \{0.4,0.8,1.2,1.9\} \\
width, depth & 256, 5 & 256, 5 & 256, 5 \\
Activation & SoftReLU & SoftReLU & SoftReLU \\
Learning rate & $1\times10^{-3}$ & $1\times10^{-3}$ & $1\times10^{-3}$ \\
sample number & 100 & 100 & 100 \\
domain mini-batch  & 128 & 128 & 128 \\
boundary mini-batch  &  &  & 64 \\
small positive number & $1\times10^{-3}$ & $1\times10^{-2}$ & $1\times10^{-2}$  \\
$r_0$ & 0.1 & 0.1 & 0.1 \\
boundary weight&  &  & 10  \\
Iterations & 50000 & 50000 & 50000 \\
\bottomrule
\end{tabularx}
\end{table}

\end{itemize}

\bibliographystyle{siamplain}

 \end{document}